\documentclass[12pt,leqno]{amsart}
\usepackage[a4paper,margin=28mm]{geometry}
\usepackage[T1]{fontenc}
\usepackage[utf8]{inputenc}
\usepackage{lmodern}
\usepackage{amsmath,amssymb,amsthm,mathtools}
\usepackage{booktabs,longtable,array,enumitem,needspace}
\usepackage{tikz-cd}
\usepackage{microtype}
\usepackage{xurl}
\usepackage[hidelinks]{hyperref}
\usepackage{xcolor}
\newif\ifshowrevisions
\showrevisionstrue

\newtheorem{theorem}{Theorem}[section]
\newtheorem{proposition}[theorem]{Proposition}
\newtheorem{lemma}[theorem]{Lemma}

\theoremstyle{definition}
\newtheorem{definition}[theorem]{Definition}
\newtheorem{remark}[theorem]{Remark}

\numberwithin{equation}{section}

\newcommand{\C}{\mathbb C}
\newcommand{\Q}{\mathbb Q}
\newcommand{\R}{\mathbb R}
\newcommand{\Z}{\mathbb Z}
\newcommand{\cO}{\mathcal O}
\newcommand{\cL}{\mathcal L}

\newcommand{\NE}{\overline{\operatorname{NE}}}
\newcommand{\Pic}{\operatorname{Pic}}
\newcommand{\G}{\operatorname{G}}
\newcommand{\Aut}{\operatorname{Aut}}
\newcommand{\Spec}{\operatorname{Spec}}

\newcommand{\Exc}{\operatorname{Exc}}

\newcommand{\red}{\mathrm{red}}

\newcommand{\wtW}{\widetilde W}
\newcommand{\wtf}{\widetilde f}

\newcommand{\pair}[2]{\langle #1,#2\rangle}
\setlist[enumerate]{label=\textup{(\roman*)},leftmargin=*,itemsep=3pt}
\hypersetup{pdftitle={Smooth minimal fourfolds not connected by smooth elementary flops}}

\title[Smooth minimal fourfolds and elementary flops]{Smooth minimal fourfolds not connected by smooth elementary flops}
\author[Osamu Fujino]{Osamu Fujino}
\author[Hiroshi Sato]{Hiroshi Sato}
\address{Department of 
Mathematics, Graduate School of Science, 
Kyoto University, Kyoto 606-8502, Japan}
\email{fujino@math.kyoto-u.ac.jp}
\address{Department of Applied Mathematics, 
Faculty of Science, Fukuoka University, 
8-19-1, Nanakuma, Jonan-ku, Fukuoka 814-0180, Japan}
\email{hirosato@fukuoka-u.ac.jp}

\date{2026/11/17, version 0.12}

\begin{document}
\begin{abstract}
We first construct two smooth projective toric fourfolds of Picard number eight that are connected by elementary flops, but cannot be connected by elementary flops through smooth complete toric varieties.
This answers a question raised by the first author in 2004.
Then, by taking suitable double covers, we construct two smooth projective minimal fourfolds of general type with Picard number nine that are connected by elementary flops, but cannot be connected by elementary flops through smooth projective varieties.
\end{abstract}
\maketitle

\section{Introduction}\label{sec:introduction}

Birational projective $\mathbb Q$-factorial terminal minimal models are connected by flops \cite[Theorem~1]{Kawamata2008Flops}. When the endpoints are smooth, it is natural to ask whether all intermediate models can also be chosen to be smooth. In dimension three, it is well known that terminal elementary flops preserve the analytic singularity type (see \cite[Theorem 2.4]{flops}). In particular, they preserve smoothness. By contrast, in dimension four or higher, an individual flop does not necessarily preserve smoothness. However, this fact alone does not answer the connectivity question. This is because a sequence passing through singular models might be replaced by a different sequence consisting only of smooth models. 

In this paper, a \emph{flop} means an elementary flop. That is, its two small contractions have relative Picard number one. This convention is essential.
Following \cite[Definition~9]{Fujino2004Fequivalence}, two smooth projective toric varieties are \emph{$F$-equivalent} if they can be connected by a finite sequence of elementary toric $K$-flops through smooth complete toric varieties. The intermediate varieties need not be projective.
We first prove the following theorem by explicitly describing toric data with the aid of generative AI. 

\begin{theorem}\label{thm:toric}
There are nonisomorphic smooth projective toric fourfolds $X_T,X_V$ of
Picard number eight and a projective $\Q$-factorial toric fourfold $X_S$
such that
\[
X_T\dashrightarrow X_S\dashrightarrow X_V
\]
is a sequence of two elementary projective $K$-flops. The only singularity
of $X_S$ is of type $\frac12(1,1,1,1)$. Neither $X_T$ nor $X_V$ admits a
nontrivial elementary toric $K$-flop 
with smooth flopped variety. 
This last statement 
also holds when the flopped variety 
is only assumed complete, rather than projective. 
In particular, $X_T$ and $X_V$ are not $F$-equivalent. 
\end{theorem}

We give all rays and maximal cones, as well as all extremal primitive relations of the smooth endpoints. Each endpoint has nine extremal rays. There are five flopping rays on $X_T$ and four on $X_V$. Every one of them has an extremal primitive relation of the form 
$v_1+v_2+v_3=w_1+2w_2$. 
Thus, all of these flops produce a singular cone of multiplicity two. 
The endpoint $X_T$ has one divisorial extremal ray, whereas $X_V$ has two. 
This fact gives a direct proof that the endpoints are not isomorphic, without assuming that the isomorphism is equivariant.

The local configuration used by the two designated flops has only six rays. It is lattice-equivalent to the configuration in \cite[Example~25]{Fujino2004Fequivalence}, which actually admits a three-step smooth factorization. Our construction adds cones that obstruct this factorization. Then we complete the local model without creating any new smooth routes. Therefore, the obstruction is not just the presence of a singular intermediate model. Rather, it is the absence of any smooth elementary modification at the very first step.

This mechanism also yields the following consequence.

\begin{proposition}\label{prop:blowup-main}
There are smooth projective toric weak Fano fourfolds
$\bar X_T,\bar X_V$ of Picard number seven, connected by six smooth
projective toric flops, and corresponding torus invariant points
$q_T\in\bar X_T$, $q_V\in\bar X_V$ in a common isomorphism open subset,
such that
\[
\operatorname{Bl}_{q_T}\bar X_T=X_T
\mbox{ and }
\operatorname{Bl}_{q_V}\bar X_V=X_V.
\]
In particular, simultaneous blow-up at corresponding points does not
preserve $F$-equivalence.
\end{proposition}

By taking suitable double covers, we obtain our main result concerning arbitrary (not necessarily toric) varieties. 

\begin{theorem}\label{thm:main}
There exist smooth projective fourfolds $M_T$ and $M_V$ over $\C$ such that:
\begin{enumerate}
\item $K_{M_T}$ and $K_{M_V}$ are nef, big and semiample, and 
$\rho(M_T)=\rho(M_V)=9$;
\item there is a sequence of four elementary projective $K$-flops from
$M_T$ to $M_V$;
\item there is no sequence of elementary $K$-flops from $M_T$ to a variety
isomorphic to $M_V$ for which every intermediate variety is smooth and
projective.
\end{enumerate}
In particular, $M_T$ and $M_V$ are of general type. 
\end{theorem}

We construct $M_T$ and $M_V$ by taking suitable double covers of $X_T$ and $X_V$, respectively. 
As mentioned above, we have explicitly determined all the extremal rays of $X_T$ and $X_V$. 
However, the behavior of the Mori cones of $M_T$ and $M_V$ is not obvious and seems very difficult to determine directly. 
Hence, Theorem~\ref{thm:main} (iii) is not a direct consequence of Theorem~\ref{thm:toric}. 
We need delicate arguments to establish Theorem~\ref{thm:main} (iii). 
The examples in Theorem~\ref{thm:main} depend on the choice of a branch divisor in a nonempty Zariski open subset of a sufficiently ample complete linear system. 
The four-step sequence in Theorem~\ref{thm:main} (ii) passes through three intermediate terminal Gorenstein varieties. 
They have one, two, and one isolated singularities of type $\frac12(1,1,1,1)$, respectively. 
We do not claim that four is the shortest possible length among all elementary flop sequences.

\subsection*{Historical background}
In 2004, the first author introduced the notion of \emph{$F$-equivalence} recalled above. 
This definition was inspired by Kawamata's paper on $D$-equivalence and $K$-equivalence (see \cite{kawamata-d-equivalence}) and Orlov's result in \cite{orlov}. In dimension three, for smooth projective toric varieties, $D$-equivalence, $K$-equivalence, and $F$-equivalence coincide. This follows easily from \cite{kawamata-d-equivalence} and \cite{orlov} (see also \cite[Theorem 15 and Remark 17]{Fujino2004Fequivalence}). We previously tried to construct an example like Theorem~\ref{thm:toric}. 
However, we could not achieve this by hand calculation alone (see \cite[Example 25]{Fujino2004Fequivalence}). 
In this paper, we finally completed this old attempt with the aid of a computer. Theorem~\ref{thm:toric} provides a counterexample to Conjectures~16 and~21 in \cite{Fujino2004Fequivalence}. It should also be noted that in 2004, \cite{BirkarEtAl2010Minimal} had not yet appeared, and a problem such as Theorem~\ref{thm:main} was completely inconceivable at that time.

\subsection*{Conventions}
We work over $\C$. 
Throughout this paper, we follow the standard notation of the minimal model theory. For details, see, for example, \cite{BirkarEtAl2010Minimal}, \cite{fujino-fundamental}, and \cite{fujino-foundations}. For toric geometry, we follow the standard textbook \cite{CoxLittleSchenck2011Toric}. 

\subsection*{Scope and verification}
The finite assertions about fans and cones are accompanied by explicit integer data and exact-arithmetic certificates. These assertions concern unimodularity, projectivity, primitive relations, wall relations, and nonnegative decompositions of curve classes. The geometric arguments about branched covers, analytic germs, and the Picard numbers of the covers are proved in the text. They are not delegated to computation. The verification code and exact data are provided in the computational supplement \cite{FujinoSato2026Verification}; their scope is described in Appendix~\ref{app:verification}. We make no claim that the Picard numbers in either theorem are minimal.

\subsection*{Organization}
Let us briefly outline the organization of this paper. 
In Section~\ref{sec:toric}, we first construct a toric variety $W$ and then explicitly give the toric data of its subdivisions $X_T$, $X_S$, and $X_V$. 
In Section~\ref{sec:extremal}, we describe the extremal and contractible primitive relations, construct the two flops, and prove 
Theorem~\ref{thm:toric}. 
In Section~\ref{sec:mechanism}, we explain some supplementary properties and prove Proposition~\ref{prop:blowup-main}. 
In Section~\ref{sec:crepant-base}, we verify the properties of the common contractions to $W$ needed for the proof of 
Theorem~\ref{thm:main}. 

Sections~\ref{sec:double-cover} to \ref{sec:automorphisms} are devoted to the proof of Theorem~\ref{thm:main}. 
In Section~\ref{sec:double-cover}, we construct suitable double covers, and in Section~\ref{sec:picard-number}, we calculate their Picard numbers. 
In Section~\ref{sec:four-step}, we prove Theorem~\ref{thm:main} (ii). 
In Section~\ref{sec:local-picard}, we discuss comparison results between the original global setting and local analytic germs. 
In Section~\ref{sec:obstruction}, we prove Theorem~\ref{thm:main} (iii) under the marking condition. 
In Section~\ref{sec:automorphisms}, we remove the marking condition assumed in Section~\ref{sec:obstruction} and finally establish 
Theorem~\ref{thm:main}. 

Finally, in Appendices~\ref{app:local-cones} to \ref{app:verification}, we provide explicit coordinate data for the fans and divisors, and explain the computer verification.

\subsection*{AI disclosure}
During the research and preparation of this paper, the authors
used ChatGPT (OpenAI, GPT-6.1) in October--November 2026 
for substantive mathematical assistance as well as manuscript
preparation.
This included exploring and developing the explicit toric examples
in Theorem~\ref{thm:toric} and
Proposition~\ref{prop:blowup-main}, and suggesting and refining
proof strategies and arguments for the branched-cover construction
and the nonconnectivity result in Theorem~\ref{thm:main}.
ChatGPT also assisted with generating and revising the computational
verification code, identifying potentially relevant literature,
and drafting and revising portions of the mathematical exposition
and English text.

The authors have critically reviewed the AI-assisted material
and checked the mathematical arguments, references, computational
procedures, and results presented in the final version.
The finite computational assertions are verified by exact-arithmetic
computations, whose scope, input data, and verification code are
described in Appendix~\ref{app:verification}.
These reproducible computations are distinct from the AI-assisted
process of discovery.
The authors take full responsibility for the correctness,
originality, and integrity of the entire work, including its
proofs, computations, and references.

\subsection*{Acknowledgments}
The first author was partially supported by JSPS KAKENHI Grant Number JP23K20787. The second author was partially supported by JSPS KAKENHI Grant Number JP24K06679.

\section{The toric construction}\label{sec:toric}
In this section, we first fix the conventions used throughout this paper. Then, we construct the common base $W$ and describe the subdivisions defining $X_T$, $X_S$, and $X_V$. 

\subsection{Conventions}
Let us start with the definition of flops. 

\begin{definition}[Flops]
Let $X$ and $X^+$ be normal $\Q$-factorial varieties. An \emph{elementary $K$-flop} is a birational map $X\dashrightarrow X^+$ satisfying the following conditions:
\begin{enumerate}
\item there exist a normal variety $P$ and projective small birational contractions $X\to P$ and $X^+\to P$,
\item the relative Picard numbers $\rho(X/P)$ and $\rho(X^+/P)$ are both one,
\item the canonical divisor $K_X$ is numerically trivial over $P$, and
\item the map $X\dashrightarrow X^+$ is the flip of some divisor on $X$.
\end{enumerate}
Note that, for the global projective sequences in this paper, $X$, $X^+$ and $P$ are all projective over $\mathbb C$. 
\end{definition}

Let $N=\Z^4$ be a lattice and let $M$ be its dual lattice. We write $\langle v_1,\ldots,v_k\rangle$ for the cone generated by the vectors $v_1,\ldots,v_k\in N$. We use the standard toric dictionary of \cite{CoxLittleSchenck2011Toric}. 
For a fan $\Sigma$, let $\G(\Sigma)$ be the set of primitive lattice generators of its one-dimensional cones. 
For a simplicial fan $\Sigma$ with convex support, a subset
$P=\{v_1,\ldots,v_h\}\subset \G(\Sigma)$ is a \emph{primitive collection} if
$\langle v_1,\ldots,v_h\rangle\notin\Sigma$, whereas every proper subset of $P$ generates a cone of $\Sigma$.
The sum $v_1+\cdots+v_h$ lies in the relative interior of a unique cone
$\langle w_1,\ldots,w_k\rangle\in\Sigma$. Its \emph{primitive relation} is
\[
v_1+\cdots+v_h=b_1w_1+\cdots+b_kw_k\quad (b_j\in\Q_{>0}).
\]
The right-hand side is zero when $k=0$. 
If $\Sigma$ is nonsingular, then all $b_j$ are integers. 

For a smooth complete toric variety $X=X_\Sigma$, let $N_1(X)_{\Z}$ denote the group of integral one-cycles modulo numerical equivalence. There is an exact sequence
\[
0\longrightarrow N_1(X)_{\Z}
\xrightarrow{\ \iota\ }\Z^{\G(\Sigma)}
\xrightarrow{\ \nu\ }N\longrightarrow0,
\]
where
\[
\iota([C])=(D_u\cdot C)_{u\in \G(\Sigma)},\quad
\nu\left((r_u)_{u\in\G(\Sigma)}\right)=\sum_{u\in \G(\Sigma)}r_u u,
\]
and $D_u$ is the torus invariant prime divisor corresponding to $u$. 
Thus, we have $N_1(X)_{\Z}\cong \ker\nu.$ 
This identifies integral numerical curve classes with integral relations among the primitive generators (see \cite[Proposition~1.3]{Casagrande2003Contractible}). 
Tensoring with $\R$ gives the corresponding description of $N_1(X)_{\R}$.
For a complete simplicial toric variety, the same description holds over $\R$ (see \cite[Proposition~6.4.1]{CoxLittleSchenck2011Toric}). 
Under this identification, 
for $r=(r_u)_{u\in\G(\Sigma)}\in \ker\nu\cong N_1(X)_{\Z}$, 
we have 
$(\sum_{u\in\G(\Sigma)}a_uD_u)\cdot r=\sum_{u\in\G(\Sigma)} a_ur_u$. 
In particular, 
$-K_X\cdot r=\sum r_u$. 
For a smooth projective toric variety, the Mori cone is generated by the wall relations. It is also generated by the primitive relations. 
For the relative quasi-projective setting, we use 
the corresponding description of relative Mori cones 
given in \cite{CoxVonRenesse2009Primitive}.

\subsection{The generators and the common base}
We take the following primitive lattice vectors
\begin{equation}\label{eq:rays}
\begin{array}{lll}
u_{1}=(-1,-1,-1,1),\ &u_{2}=(-1,0,1,1),\ 
&u_{3}=(0,-1,-1,1),\\
u_{4}=(0,0,-1,1),\ 
&u_{5}=(0,0,0,1),\ &u_{6}=(0,1,0,1),\\
u_{7}=(1,-1,-1,1),\ &u_{8}=(1,1,0,1),\ 
&u_{9}=(0,0,-1,0),\\ 
u_{10}=(0,0,0,-1),\ 
&u_{11}=(1,0,-1,1),\ &u_{12}=(3,0,-3,4)
\end{array}
\end{equation}
in $N$. 
The first eight vectors lie in the hyperplane where the last coordinate is one. We set 
$\sigma:=\langle u_{1},\ldots,u_{8}\rangle$. 
Since 
$2u_3=u_1+u_7$ and $5u_5=2u_2+u_4+u_7+u_8$, 
we have $\sigma=\langle u_1,u_2,u_4,u_6,u_7,u_8\rangle$. 
The six displayed generators are precisely 
the primitive ray generators of $\sigma$. 
We first construct the common global fan $\Delta_W$.
Its maximal cones are $\sigma$ and the following eighteen cones:
\begin{equation}\label{eq:base-maximal-cones}
\begin{array}{lll}
\langle u_{1},u_{2},u_{7},u_{10}\rangle,&
\langle u_{1},u_{4},u_{7},u_{9}\rangle,&
\langle u_{1},u_{7},u_{9},u_{10}\rangle,\\
\langle u_{1},u_{4},u_{6},u_{9}\rangle,&
\langle u_{2},u_{6},u_{8},u_{10}\rangle,&
\langle u_{2},u_{7},u_{8},u_{10}\rangle,\\
\langle u_{4},u_{6},u_{8},u_{9}\rangle, &
\langle u_{4},u_{7},u_{8},u_{12}\rangle,&
\langle u_{4},u_{7},u_{9},u_{11}\rangle,\\
\langle u_{4},u_{7},u_{11},u_{12}\rangle,&
\langle u_{4},u_{8},u_{9},u_{11}\rangle,&
\langle u_{4},u_{8},u_{11},u_{12}\rangle,\\
\langle u_{6},u_{8},u_{9},u_{10}\rangle,&
\langle u_{7},u_{8},u_{10},u_{11}\rangle,&
\langle u_{7},u_{8},u_{11},u_{12}\rangle,\\
\langle u_{7},u_{9},u_{10},u_{11}\rangle,&
\langle u_{8},u_{9},u_{10},u_{11}\rangle,&
\langle u_{1},u_{2},u_{6},u_{9},u_{10}\rangle.
\end{array}
\end{equation}
We include all faces of these nineteen cones. We set $W=X_{\Delta_W}$ and denote by $p\in W$ the torus invariant point corresponding to $\sigma$.
The primitive generators of the rays of $\Delta_W$ are the vectors in \eqref{eq:rays} other than $u_{3}$ and $u_{5}$.
The completeness and projectivity assertions are verified in Proposition~\ref{prop:fan-check} below.

\subsection{The three subdivisions}
We first make two subdivisions that are common to all three models 
$X_T$, $X_S$ and $X_V$, 
and then refine the cones supported on $\sigma$ 
in three different ways in order to construct the fans 
associated to $X_T$, $X_S$ and $X_V$. 
At each stage, we include all faces of the specified maximal cones.

\smallskip
\noindent\emph{First common subdivision.}
We replace $\sigma$ by the following five maximal cones:
\begin{equation}\label{eq:common-local-subdivision}
\begin{aligned}
\tau_1&=\langle u_{1},u_{2},u_{3},u_{4},u_{5},u_{6}\rangle,\ 
\tau_2=\langle u_{2},u_{3},u_{5},u_{7}\rangle,\ 
\tau_3=\langle u_{2},u_{5},u_{6},u_{8}\rangle,\\
\tau_4&=\langle u_{2},u_{5},u_{7},u_{8}\rangle,\ 
\tau_5=\langle u_{3},u_{4},u_{5},u_{6},u_{7},u_{8}\rangle.
\end{aligned}
\end{equation}
This introduces the rays generated by $u_3$ and $u_5$.
Since $2u_3=u_1+u_7$, the subdivision of the boundary of $\sigma$ inserts $u_3$ into the two-dimensional cone $\langle u_1,u_7\rangle$.
Accordingly, we simultaneously subdivide the three exterior maximal cones containing this two-dimensional cone as follows:
\[
\begin{aligned}
\langle u_1,u_2,u_7,u_{10}\rangle
&\longmapsto
\{\langle u_1,u_2,u_3,u_{10}\rangle,\langle u_2,u_3,u_7,u_{10}\rangle\},\\
\langle u_1,u_4,u_7,u_9\rangle
&\longmapsto
\{\langle u_1,u_3,u_4,u_9\rangle,\langle u_3,u_4,u_7,u_9\rangle\},\\
\langle u_1,u_7,u_9,u_{10}\rangle
&\longmapsto
\{\langle u_1,u_3,u_9,u_{10}\rangle,\langle u_3,u_7,u_9,u_{10}\rangle\}.
\end{aligned}
\]
All other maximal cones of $\Delta_W$ remain unchanged.
These subdivisions agree on common faces and define a fan $\Delta^{(1)}$ with twenty-six maximal cones.
This is a simultaneous subdivision of $\sigma$ and the three adjacent exterior cones, not a sequence of star subdivisions at $u_3$ and $u_5$.

\smallskip
\noindent\emph{Second common subdivision.}
The cone $\kappa=\langle u_1,u_2,u_6,u_9,u_{10}\rangle$ is still a maximal cone of $\Delta^{(1)}$.
We replace it by
\[
\langle u_1,u_2,u_6,u_9\rangle,\quad
\langle u_1,u_2,u_9,u_{10}\rangle,\quad
\langle u_2,u_6,u_9,u_{10}\rangle.
\]
This subdivision, corresponding to $u_1+u_6+u_{10}=u_2+2u_9$, introduces no new ray and leaves every proper face of $\kappa$ unchanged.
Keeping all other maximal cones gives a common fan $\Delta^{(2)}$ with twenty-eight maximal cones.
Its five maximal cones contained in $\sigma$ are still $\tau_1,\ldots,\tau_5$.

\smallskip
\noindent\emph{The three final subdivisions.}
We now refine $\tau_1$ and $\tau_5$ in \eqref{eq:common-local-subdivision} in three ways. To begin with, 
the five cones
\[
\begin{gathered}
\langle u_1,u_2,u_3,u_5\rangle,\quad
\langle u_1,u_2,u_4,u_5\rangle,\quad
\langle u_1,u_2,u_4,u_6\rangle,\\
\langle u_1,u_3,u_4,u_5\rangle,\quad
\langle u_2,u_4,u_5,u_6\rangle
\end{gathered}
\]
give the same subdivision of $\tau_1$ in all three cases. 

Next, 
let $T,S,V$ denote the following sets of maximal cones of subdivisions of $\tau_5$:
\begin{equation}\label{eq:local-fans}
T=\left\{\langle u_{3},u_{4},u_{5},u_{8}\rangle,
\langle u_{3},u_{4},u_{7},u_{8}\rangle,\langle u_{3},u_{5},u_{7},u_{8}\rangle,
\langle u_{4},u_{5},u_{6},u_{8}\rangle \right\}.
\end{equation}
\begin{equation}\label{eq:local-S}
S=\left\{\langle u_{3},u_{4},u_{5},u_{7}\rangle,
\langle u_{4},u_{5},u_{6},u_{8}\rangle,\langle u_{4},u_{5},u_{7},u_{8}\rangle\right\}.
\end{equation}
\begin{equation}\label{eq:local-V}
V=\left\{\langle u_{3},u_{4},u_{5},u_{7}\rangle,
\langle u_{4},u_{5},u_{6},u_{7}\rangle,\langle u_{4},u_{6},u_{7},u_{8}\rangle,\langle u_{5},u_{6},u_{7},u_{8}\rangle\right\}.
\end{equation}
These refinements agree on common faces and do not change the subdivision of the boundary of $\sigma$. 
Outside $\sigma$, we keep the subdivision already made in $\Delta^{(2)}$.
Finally, we obtain the three fans $\Sigma_T$, $\Sigma_S$ and $\Sigma_V$ corresponding to 
these subdivisions $T$, $S$ and $V$, respectively. 
Their corresponding toric varieties are $X_T,X_S,X_V$. 
On $X_T,X_S,X_V$, we denote by $D_i=D_{u_i}$ the torus
invariant prime divisor corresponding to $u_i$. 
$\Sigma_T$, $\Sigma_S$ and $\Sigma_V$ refine $\Delta^{(2)}$, which refines $\Delta^{(1)}$, which in turn refines $\Delta_W$.
In particular, all three fans subdivide $\Delta_W$, and 
$\G(\Sigma_T)=\G(\Sigma_S)=\G(\Sigma_V)=\{u_1,\ldots,u_{12}\}$. 
These assertions can be checked directly from the cones listed above. 
Put $W_i=X_{\Delta^{(i)}}$ for $i=1,2$. The identity of $N$ induces proper birational toric morphisms
\begin{equation}\label{eq:common-morphism-chain}
X_\lambda\longrightarrow W_2\longrightarrow W_1\longrightarrow W
\quad(\lambda=T,S,V).
\end{equation}
We denote their composites by $f_\lambda:X_\lambda\to W$. 
For $\lambda=T,S,V$, we set 
$\Sigma_\lambda|_\sigma
:=\{\omega\in\Sigma_\lambda\mid \omega\subseteq\sigma\}$. 
This is the subfan of $\Sigma_\lambda$ supported on $\sigma$. 

\begin{proposition}\label{prop:fan-check}
The fans $\Delta_W,\Sigma_T,\Sigma_S,\Sigma_V$ are projective. 
The fans $\Sigma_T$ and $\Sigma_V$ have 
$35$ maximal cones, all of which are nonsingular. The fan $\Sigma_S$ has
$34$ maximal cones. Its only singular maximal cone is
$\langle u_4,u_5,u_7,u_8\rangle$, which has multiplicity two.
In particular, $X_T$ and $X_V$ are smooth and their Picard numbers are eight.
\end{proposition}
\begin{proof}
The assertions on completeness, projectivity, nonsingularity, and multiplicities
are verified by exact computer calculations from the displayed fan data.
See Appendix~\ref{app:verification} and the ancillary verification files.
Finally, $\rho(X_T)=\rho(X_V)=12-4=8$.
\end{proof}

\begin{remark}
Let us describe the morphisms in
\eqref{eq:common-morphism-chain}.
For $\lambda=T,S,V$, the morphism
$X_\lambda\to W_2$ is small and has relative
Picard number three.
The morphism $W_2\to W_1$ is a small
contraction of relative Picard number one.
It contracts the invariant surface
$V(\langle u_2,u_9\rangle)$ to the point
corresponding to $\kappa$.
Finally, $W_1\to W$ has relative Picard
number one, but contracts two invariant
prime divisors. The divisors corresponding
to $u_3$ and $u_5$ are mapped onto
$V(\langle u_1,u_7\rangle)$ and $p$,
respectively.
\end{remark}

\section{Extremal primitive relations and the obstruction}\label{sec:extremal}
\subsection{The complete Mori cones}\label{completemori3.1}
We recall the description of contractions 
associated with extremal primitive relations. 
See \cite[Corollary~2.4]{Casagrande2003Contractible} for details. 
Let $X=X_\Sigma$ be a smooth projective toric fourfold, and let an extremal ray be represented by the primitive relation
\[
r:\ v_1+\cdots+v_h=b_1w_1+\cdots+b_kw_k
\quad (b_j\in\Z_{>0}).
\]
Its degree is $h-\sum_jb_j=(-K_X)\cdot r$.
If $k=0$, the contraction is a $\mathbb P^{h-1}$-bundle over a smooth variety
of dimension $5-h$.
If $k>0$, the contraction $\varphi:X\to Y$ is birational, and
\[
\Exc(\varphi)=V(\langle w_1,\ldots,w_k\rangle),\quad
\dim\Exc(\varphi)=4-k,\quad
\dim\varphi(\Exc(\varphi))=5-h-k.
\]
In particular, $k=1$ gives a divisorial contraction, whereas $k\ge2$ gives a
small contraction. The birational contraction target $Y$ is smooth if and only if
$k=1$ and $b_1=1$. In this case, $\varphi$ is the blow-up of a smooth torus
invariant center. 
For a small contraction, the variety $X^+$ obtained by the flip, anti-flip,
or flop is smooth if and only if $b_1=\cdots=b_k=1$. This follows from the
fan description and the corresponding cone multiplicities. 
A small extremal contraction is flopping precisely when its degree is zero.

We also use contractible classes in the sense of
\cite[Definition~2.3]{Casagrande2003Contractible}.
A primitive relation $r$ is \emph{contractible} if some irreducible curve
has numerical class in $\R_{>0}r$ and there is a toric morphism
$\varphi_r:X\to Y$ with connected fibers such that 
an irreducible curve $C\subset X$ 
is contracted by $\varphi_r$ if and only if 
$[C]\in\R_{\ge0}r$. 
The contraction target is not required to be projective.
Every extremal primitive relation is contractible. 
For the fans $\Sigma_T$ and $\Sigma_V$, the converse also holds, 
as we shall verify below. 

The extremal primitive relations of $\Sigma_T$ are
$\alpha_1,\ldots,\alpha_6,\beta_1,\beta_2,\beta_3$, and those of $\Sigma_V$
are $\alpha_1,\ldots,\alpha_6,\gamma_1,\gamma_2,\gamma_3$, where
\begin{equation}\label{eq:common-extremal-relations}
\begin{array}{ll}
\alpha_1:\ u_9+u_{12}=u_4+3u_{11},
 &\alpha_2:\ u_1+u_5+u_6=u_2+2u_4,\\
\alpha_3:\ u_1+u_6+u_{10}=u_2+2u_9,
 &\alpha_4:\ u_2+u_3+u_4=u_1+2u_5,\\
\alpha_5:\ u_2+u_4+u_9=u_1+u_6,
&\alpha_6:\ u_4+u_7+u_8+u_{11}=u_{12},
\end{array}
\end{equation}
and
\begin{equation}\label{eq:other-extremal-relations}
\begin{array}{ll}
\beta_1:\ u_3+u_{12}=2u_4+2u_7+u_8,
& \beta_2:\ u_2+u_4+u_8=2u_5+u_6,\\
\beta_3:\ u_4+u_5+u_7=2u_3+u_8, & \\
\gamma_1:\ u_6+u_{12}=2u_4+u_7+2u_8,
& \gamma_2:\ u_1+u_7=2u_3,\\
\gamma_3:\ u_4+u_5+u_8=2u_6+u_7. &
\end{array}
\end{equation}

\begin{proposition}\label{prop:mori-cones}
We consider the complete Mori cones inside the common vector space of relations among the twelve rays. Then we have
\[
\NE(X_T)=
\langle \alpha_1,\ldots,\alpha_6,\beta_1,\beta_2,\beta_3\rangle
\mbox{ and }
\NE(X_V)=\langle \alpha_1,\ldots,\alpha_6,\gamma_1,\gamma_2,\gamma_3\rangle.
\]
All displayed generators are extremal, and both cones have dimension eight.
For each endpoint, these nine primitive relations are also exactly the contractible primitive relations. 
The flopping extremal rays are exactly 
$\alpha_2$, $\alpha_3$, $\alpha_4$, $\beta_2$, $\beta_3$ on $X_T$, 
and 
$\alpha_2$, $\alpha_3$, $\alpha_4$, $\gamma_3$ 
on $X_V$. 
\end{proposition}
\begin{proof}
The complete lists and their extremality are verified by exact computer
calculations. Contractibility is checked using
\cite[Proposition~3.4]{Casagrande2003Contractible}.
Appendix~\ref{app:primitive-relations} gives all primitive relations, grouped
by these properties. The verification data are described in
Appendix~\ref{app:verification}.
The preceding descriptions of extremal contractions give the asserted
flopping rays.
\end{proof}

\begin{remark}
$u_{1}+u_{7}=2u_{3}$ is primitive on both endpoints. However, 
its class on $X_T$ is given by 
$\gamma_2=\alpha_2+\beta_2+\beta_3$. 
Therefore, it is nonextremal on $X_T$ and becomes extremal on $X_V$.
The designated two-step factorization uses the oriented rays $\beta_3$ and
$-\gamma_3$.
\end{remark}

For completeness, we note that the intermediate variety $X_S$ has $68$ walls, $31$ primitive
collections, and $11$ extremal rays. We set 
$\delta:\ u_{5}+2u_{12}=3u_{4}+3u_{7}+3u_{8}$. 
Then we have
\[
\NE(X_S)=\langle \alpha_1,\ldots,\alpha_6,\beta_2,\gamma_2,
-\beta_3,-\gamma_3,\delta\rangle.
\]
Since $X_S$ is singular, the displayed generators are 
primitive integral vectors in the relation lattice. 
They need not be the numerical classes of individual 
torus invariant curves. 
The same exact wall certificates verify this list.

\subsection{The two flops}\label{subsec:two-flops}
Let $P_{TS}$ and $P_{SV}$ be the targets of the extremal contractions
associated with $\beta_3$ on $X_T$ and $\gamma_3$ on $X_V$, respectively.
Both contractions are small and $K$-trivial, and their targets are projective. 
Comparing the maximal cones in \eqref{eq:local-fans}, 
\eqref{eq:local-S} and \eqref{eq:local-V} 
shows that the flopped variety in either case is $X_S$. 
Namely, the flop associated with 
the primitive relation $\beta_3$ on $X_T$ replaces the cones 
\[
\langle u_{3},u_{4},u_{5},u_{8}\rangle,\ 
\langle u_{3},u_{4},u_{7},u_{8}\rangle,\ 
\langle u_{3},u_{5},u_{7},u_{8}\rangle\in T
\]
by the cones 
\[
\langle u_{3},u_{4},u_{5},u_{7}\rangle,\ 
\langle u_{4},u_{5},u_{7},u_{8}\rangle \in S,
\]
while the flop associated with 
the primitive relation $\gamma_3$ on $X_V$ replaces the cones 
\[
\langle u_{4},u_{5},u_{6},u_{7}\rangle,\ 
\langle u_{4},u_{6},u_{7},u_{8}\rangle,\ 
\langle u_{5},u_{6},u_{7},u_{8}\rangle \in V
\]
by the cones 
\[
\langle u_{4},u_{5},u_{6},u_{8}\rangle,\ 
\langle u_{4},u_{5},u_{7},u_{8}\rangle\in S. 
\]
Reversing the flop from $X_V$, we obtain 
\[
X_T\dashrightarrow X_S\dashrightarrow X_V.
\]
Both contraction targets carry morphisms to $W$: their fans are common
coarsenings of the adjacent fans and still subdivide $\Delta_W$.
Thus all the morphisms in
\[
X_T\longrightarrow P_{TS}\longleftarrow X_S,
\quad
X_S\longrightarrow P_{SV}\longleftarrow X_V
\]
are over $W$. 

On $X_S$, the nontrivial lattice class of the cone $\langle u_{4},u_{5},u_{7},u_{8}\rangle\in\Sigma_S$ is represented by the element
\[
\frac{u_{4}+u_{5}+u_{7}+u_{8}}{2}=u_{3}+u_{8}\in N.
\]
In this expression, the four coefficients on the left-hand side are all one half, which means its lattice index is two.
Thus, its corresponding affine chart has a quotient singularity of type $\frac12(1,1,1,1)$.
All proper faces of this cone are nonsingular, so this singularity is isolated. The nonidentity group element has age two. Consequently, the singularity is terminal and Gorenstein. This completes the proof of the two-step singular factorization in Theorem~\ref{thm:toric}.

\subsection{Every flopping ray leads to a singular model}
\label{floppingsing3.3}
By \eqref{eq:common-extremal-relations} and 
\eqref{eq:other-extremal-relations}, all flopping
extremal primitive relations of $X_T$ and $X_V$
are of the form 
$v_1+v_2+v_3=w_1+2w_2$. 
As explained in Section~\ref{completemori3.1}, the corresponding
flops necessarily produce singular varieties. 
Thus, neither $X_T$ nor $X_V$ admits an 
elementary projective toric $K$-flop whose 
flopped variety is smooth.

We also notice that $X_T$ has exactly one
divisorial extremal ray, while $X_V$ has exactly two. Any abstract isomorphism
between them would identify their Mori cones and preserve the dimensions of the exceptional
loci of their extremal contractions. Therefore, we conclude that 
$X_T\not\simeq X_V$ 
as abstract varieties. This conclusion relies entirely on the complete extremal-ray calculation. It does not assume that an isomorphism respects the chosen tori.

\subsection{The relative Mori cones at the distinguished point}\label{subsec:relative}
Let $Z_T$ and $Z_V$ be the toric varieties associated with 
the fans $\Sigma_T|_\sigma$ and $\Sigma_V|_\sigma$, respectively. 
They come with projective birational morphisms to $U_\sigma$. 
The local calculation needed for the double covers is
\begin{equation}\label{eq:relative-cones}
\NE(Z_T/U_\sigma)=\langle \alpha_2,\alpha_4,\beta_2,\beta_3\rangle,
\quad 
\NE(Z_V/U_\sigma)=\langle \alpha_2,\alpha_4,\gamma_2,\gamma_3\rangle.
\end{equation}
Each of these cones is simplicial of dimension four. To verify
\eqref{eq:relative-cones}, we enumerate the $19$ interior walls of each local
fan. They give nine distinct relation vectors. All of them are nonnegative
integral combinations of the four displayed generators. These four generators are linearly independent and they themselves occur as wall classes.
Appendix~\ref{app:local-cones} gives these decompositions explicitly.
Only the interior walls produce complete invariant curves in the local model. These curves generate the relative Mori cone (see, for example, \cite[Corollary 6.3.21]{CoxLittleSchenck2011Toric}).

Every small extremal ray of $\NE(Z_T/U_\sigma)$
and $\NE(Z_V/U_\sigma)$ corresponds to a primitive relation
of the form 
$v_1+v_2+v_3=w_1+2w_2$. 
As explained in Section~\ref{completemori3.1}, the corresponding flop
necessarily produces a singular variety.
Hence, neither local model has a nontrivial
smooth elementary flop over $U_\sigma$. This is precisely the local statement used in the proof of Proposition~\ref{prop:marked-nonconnectivity}, in Section~\ref{subsec:numerical-localization}, to exclude non-toric smooth projective factorizations.

\subsection{Allowing nonprojective toric intermediates}
By \cite[Theorem~2.2 and Lemma~2.6]{Casagrande2003Contractible}, every
nontrivial elementary toric contraction from a smooth complete toric
variety is the contraction of a contractible primitive relation, even when
the contraction target is not projective. Proposition~\ref{prop:mori-cones} shows that,
on either endpoint, the contractible primitive relations are exactly the
extremal ones. Consequently, every possible first elementary toric
$K$-flop is among those already excluded in 
Section~\ref{floppingsing3.3}. 
This proves the nonprojective assertion and completes the proof of
Theorem~\ref{thm:toric}.

\section{\texorpdfstring{$F$}{F}-equivalence and point blow-ups}\label{sec:mechanism}
Our example is designed to obstruct a smooth factorization, rather than
merely to exhibit a singular intermediate model. In this section, we explain both the
local and global aspects of this distinction. These observations are not used in the subsequent construction or proof of Theorem~\ref{thm:main}.

\subsection{A six-ray configuration that does admit smooth flops}
Only the six rays $u_{3},u_{4},u_{5},u_{6},u_{7},u_{8}$ occur in the two designated
relations $\beta_3$ and $\gamma_3$. 
The two subdivisions of $\tau_5$ with maximal cones $T$ and $V$ 
give toric varieties connected by the following smooth 
elementary toric flops: 
\begin{equation}\label{eq:six-ray-path}
\begin{aligned}
(1)&\quad u_3+u_6=u_4+u_5,\\
(2)&\quad u_6+u_7=u_3+u_8,\\
(3)&\quad u_4+u_5=u_3+u_6.
\end{aligned}
\end{equation}
Each relation is written with the primitive collection on the left in the
fan immediately before the corresponding flop.
The complete lists of maximal cones along this path are as follows:
\[
\begin{array}{cl}
&\left\{\langle u_{3},u_{4},u_{5},u_{8}\rangle,
\langle u_{3},u_{4},u_{7},u_{8}\rangle,
\langle u_{3},u_{5},u_{7},u_{8}\rangle,
\langle u_{4},u_{5},u_{6},u_{8}\rangle\right\}=T\\
\overset{(1)}\dashrightarrow&\left\{\langle u_{3},u_{4},u_{6},u_{8}\rangle,
\langle u_{3},u_{4},u_{7},u_{8}\rangle,
\langle u_{3},u_{5},u_{6},u_{8}\rangle,
\langle u_{3},u_{5},u_{7},u_{8}\rangle\right\}\\
\overset{(2)}\dashrightarrow&
\left\{\langle u_{3},u_{4},u_{6},u_{7}\rangle,
\langle u_{3},u_{5},u_{6},u_{7}\rangle,
\langle u_{4},u_{6},u_{7},u_{8}\rangle,
\langle u_{5},u_{6},u_{7},u_{8}\rangle\right\}\\
\overset{(3)}\dashrightarrow&
\left\{\langle u_{3},u_{4},u_{5},u_{7}\rangle,
\langle u_{4},u_{5},u_{6},u_{7}\rangle,
\langle u_{4},u_{6},u_{7},u_{8}\rangle,
\langle u_{5},u_{6},u_{7},u_{8}\rangle\right\}=V.
\end{array}
\]
The indicated cone replacements define toric flops, and all the cones are nonsingular. Furthermore, the models and their contraction targets are 
projective over the affine toric variety associated with $\tau_5$. 
These finite checks are included in 
the verification data. 

This is precisely the six-ray configuration presented in
\cite[Example~25]{Fujino2004Fequivalence}. In the notation of those notes,
the correspondence is given by
\[
(x_1,x_2,x_3,x_4,x_6,x_8)
\longmapsto(u_{4},u_{3},u_{5},u_{6},u_{7},u_{8}).
\]
It is induced by the lattice automorphism
\[
(x,y,z,t)\longmapsto (z,-x+y,y-t,t),
\]
which has determinant $-1$. Thus, the central two-step singular path is
not, by itself, the source of the obstruction. 

In the eight-ray fan $\Sigma_T|_\sigma$, the additional 
maximal cone $\langle u_1,u_2,u_4,u_5\rangle$ contains the face 
$\langle u_4,u_5\rangle$ but contains neither $u_3$ nor $u_6$. 
This cone blocks the first exchange of
\eqref{eq:six-ray-path}. The relative Mori cone calculation \eqref{eq:relative-cones} also excludes all other possible first smooth steps.
This explains the role of the rays that do not appear in the two designated flop
relations.

\subsection{Removing the common point blow-up}

Recall that the extremal primitive relation
\begin{equation}\label{eq:point-blowup}
\alpha_6:\ u_{4}+u_{7}+u_{8}+u_{11}=u_{12}
\end{equation}
is common to $\Sigma_T$ and $\Sigma_V$.
For $\lambda=T,V$, the corresponding extremal contraction 
$X_\lambda\to\bar X_\lambda$ 
is the blow-up at a nonsingular torus invariant point 
corresponding to the cone $\langle u_4,u_7,u_8,u_{11}\rangle$. 
We denote the fan of $\bar X_\lambda$ by
$\bar\Sigma_\lambda$.

\begin{proposition}\label{prop:weak-fano-path}
The varieties $\bar X_T$ and $\bar X_V$ are smooth projective toric weak Fano
fourfolds with Picard number seven. They are connected by a sequence of six elementary
flops through smooth projective toric weak Fano varieties.
\end{proposition}
\begin{proof}
The varieties $\bar X_T$ and $\bar X_V$ are smooth 
and projective by the description of the extremal 
contraction associated with $\alpha_6$. 
Their Picard numbers are seven. 
All their primitive relations, 
listed in Appendix~\ref{app:primitive-relations}, 
have nonnegative degree, so both varieties are weak Fano.

Starting with $\bar\Sigma_T$, we perform the flops associated successively
with the following primitive relations:
\begin{equation}\label{eq:weak-fano-path}
\begin{array}{ll}
(1)\ u_3+u_{11}=u_4+u_7, & (2)\ u_5+u_{11}=u_3+u_8,\\
(3)\ u_4+u_7=u_3+u_{11}, & (4)\ u_6+u_{11}=u_4+u_8,\\
(5)\ u_6+u_7=u_5+u_{11}, & (6)\ u_4+u_8=u_6+u_{11}.
\end{array}
\end{equation}
Appendix~\ref{app:weak-fano} lists the four maximal cones
removed and the four maximal cones inserted at each step.
Exact computer calculations verify that all seven fans
are complete, projective, and nonsingular, that the
corresponding primitive relations are extremal, and
that every primitive relation at every stage has
nonnegative degree. 
Since each transformation is small and crepant,
all six steps are elementary projective $K$-flops
through smooth projective toric weak Fano fourfolds.
The final fan is $\bar\Sigma_V$.
\end{proof}

The birational map $\bar X_T\dashrightarrow\bar X_V$ defined by the original two
singular flops acts as the identity on the affine chart corresponding to $\langle u_4,u_7,u_8,u_{11}\rangle$. 
The two torus invariant points therefore represent corresponding points in a common
isomorphic open subset. Blowing them up yields $X_T$ and $X_V$, which are not
$F$-equivalent by Theorem~\ref{thm:toric}. This proves
Proposition~\ref{prop:blowup-main}. In particular, no smooth elementary
factorization of this map can avoid the chosen point throughout its entire sequence; if it did, such a
factorization would lift to a smooth factorization between the blow-ups.

\begin{remark}\label{rem:minimality}
The fourfolds $X_T$ and $X_V$ are not weak Fano. For example, we have
$\deg\beta_1=-3$ and $\deg\gamma_1=-3$. 
By Proposition~\ref{prop:mori-cones}, the only extremal contraction of either endpoint that is the inverse of a blow-up along a smooth invariant center is the contraction associated with $\alpha_6$, namely \eqref{eq:point-blowup}. Indeed, the only other divisorial extremal ray is $\gamma_2$ on $X_V$, whose 
contraction target is singular. 
Removing the blow-up \eqref{eq:point-blowup} restores smooth connectivity. This observation does not prove that Picard number eight is the minimum among all examples of obstructed smooth toric flops.
\end{remark}

\section{The common crepant contractions}\label{sec:crepant-base}
We now verify the properties of the common base $W$, constructed in Section~\ref{sec:toric}, that are needed for the double covers.

\begin{proposition}\label{prop:common-base}
The variety $W$ is a projective Gorenstein canonical toric fourfold.
For each ray of $\Delta_W$ generated by $u_i$, let $D_i^W$ be the corresponding invariant prime divisor on $W$.
The divisor
\[
A=6D_9^W+12D_{10}^W+2D_{11}^W+D_{12}^W
\]
is an ample Cartier divisor, and there are projective birational toric morphisms
\[
f_\lambda:X_\lambda\longrightarrow W\quad(\lambda=T,S,V)
\]
such that
\begin{equation}\label{eq:crepant-base}
K_{X_\lambda}=f_\lambda^*K_W.
\end{equation}
Moreover, $K_W+4A$ is an ample Cartier divisor. The cone $\sigma$ is a maximal cone
of $\Delta_W$.
\end{proposition}
\begin{proof}
Exact computer calculations verify that $A$ and $K_W$ are Cartier, that
$A$ and $K_W+4A$ are ample, and that the subdivisions give
$K_{X_\lambda}=f_\lambda^*K_W$ for $\lambda=T,S,V$.
See Appendices~\ref{app:base} and~\ref{app:verification} and the ancillary files.
The proper birational morphisms $f_\lambda$ constructed in
\eqref{eq:common-morphism-chain} are projective since their sources are
projective. Normal toric varieties are Cohen--Macaulay, so the Cartier
property of $K_W$ implies that $W$ is Gorenstein. The smooth crepant
resolution $f_T$ implies that $W$ has canonical singularities.
Finally, $\sigma$ is a maximal cone of $\Delta_W$ by construction.
\end{proof}

Recall that $p\in W$ is the torus invariant point corresponding to $\sigma$. Its invariant affine neighborhood is 
$U_\sigma=\Spec\C[\sigma^\vee\cap M]$. 
The inverse images $Z_T=f_T^{-1}(U_\sigma)$ and 
$Z_V=f_V^{-1}(U_\sigma)$ are exactly the toric varieties 
associated with the local fans $\Sigma_T|_\sigma$ and 
$\Sigma_V|_\sigma$, respectively. 
The last-coordinate function takes the value one on every generator in $\sigma$. In particular, the local resolutions are crepant.

\begin{lemma}\label{lem:divisor-images}
For $\lambda=T,V$, the divisor $D_{5}$ is the unique torus-invariant prime divisor on $X_\lambda$ whose image under $f_\lambda$ is a point. Moreover, $f_\lambda(D_{5})=p$.
\end{lemma}
\begin{proof}
The generator $u_{5}=(2u_2+u_4+u_7+u_8)/5$ 
lies in the interior of $\sigma$. 
The generator $u_{3}=(u_{1}+u_{7})/2$ lies in the relative interior of the two-dimensional cone $\langle u_{1},u_{7}\rangle$ of $\Delta_W$. All other $u_i$ generate rays of $\Delta_W$. The orbit-cone correspondence therefore gives the assertion.
\end{proof}
\color{black}

\section{Branched double covers and minimal models}\label{sec:double-cover}
In this section, we construct suitable branched double covers.

\subsection{The branch conditions}
We take a sufficiently large and divisible positive integer $m$ such that $K_W+mA$ is very ample. We set 
$\cL=\cO_W(mA)$ and $L_\lambda=mf_\lambda^*A$. 
We first require a divisor $B\in|2mA|$ to satisfy the following conditions: 
\begin{enumerate}
\item $B$ is nonzero, reduced, and smooth on the smooth locus of $W$, 
\item $p\notin B$, and 
\item $B_T=f_T^*B$ and $B_V=f_V^*B$ are smooth.
\end{enumerate}
We call these the \emph{basic branch conditions}. These conditions hold on a nonempty Zariski open subset of $|2mA|$. 

For the assertion about abstract endpoints, we will also impose the condition 
$\Aut(W,B)=\{1\}$. 
The existence of a nonempty Zariski open subset of $|2mA|$ with this property is proved in Section~\ref{sec:automorphisms}. 

\subsection{The double cover of the singular base}
Let $s_B\in H^0(W,\cL^{\otimes2})$ be a section defining $B$. We consider the $\mathcal O_W$-algebra $\cO_W\oplus\cL^{-1}$ whose algebra structure is given by the multiplication $\cL^{-2}\xrightarrow{s_B}\cO_W$. We set
\begin{equation}\label{eq:base-cover}
\pi:\wtW=\Spec_W(\cO_W\oplus\cL^{-1})\longrightarrow W.
\end{equation}
This is a finite flat double cover locally defined by $z^2=s_B$. 
By construction, $\wtW$ is a normal Gorenstein variety. 

We put 
$M_\lambda:=X_\lambda\times_W\wtW$ $(\lambda=T,V)$. 
We obtain the following Cartesian squares:
\[
\begin{tikzcd}[column sep=large,row sep=large]
M_\lambda\arrow[r,"\pi_\lambda"]\arrow[d,"\wtf_\lambda"']&
X_\lambda\arrow[d,"f_\lambda"]\\
\wtW\arrow[r,"\pi"']&W.
\end{tikzcd}
\]
Each $\pi_\lambda$ is a double cover branched over $B_\lambda=f^*_\lambda B$. 
Thus, $M_T$ and $M_V$ are smooth integral projective fourfolds, and the maps $\wtf_\lambda$ are projective and birational.

\subsection{Canonical divisors and the canonical model}
By construction, we have 
\begin{equation}\label{eq:double-cover-adjunction}
K_{\wtW}=\pi^*\left(K_W+\frac{1}{2}B\right),\quad
K_{M_\lambda}=\pi_\lambda^*\left(K_{X_\lambda}+\frac{1}{2}B_\lambda\right).
\end{equation}
By combining \eqref{eq:double-cover-adjunction} with \eqref{eq:crepant-base}, we obtain 
\begin{equation}\label{eq:canonical-cover}
K_{M_\lambda}=\wtf_\lambda^*\pi^*\left(K_W+\frac{1}{2}B\right)
=\wtf_\lambda^*K_{\wtW}.
\end{equation}
Since the divisor $K_W+\frac{1}{2}B$ is ample and $\pi$ is finite, $K_{\wtW}$ is ample. Consequently, $K_{M_\lambda}$ is big and semiample.
Therefore, $M_T$ and $M_V$ are minimal models of general type without boundary divisors.

Note that the variety $\wtW$ has canonical singularities since it has a smooth crepant resolution. By the projection formula, we have 
$H^0(M_\lambda,rK_{M_\lambda})=H^0(\wtW,rK_{\wtW})$ for every $r\geq0$. 
This means that $\wtW$ is the common canonical model. A finite flat base change preserves the codimension of the exceptional loci. Thus, the induced birational map $M_T\dashrightarrow M_V$ is small.

\begin{remark}
The divisor $f_\lambda^*A$ is nef and big, but not ample. It has degree zero on the curves where we retain the local obstruction.
\end{remark}

\section{The Picard number of the covers}\label{sec:picard-number}
We compute the Picard number of every cover satisfying the basic branch
conditions. We write $X=X_T$ or $X_V$, $f=f_\lambda$, $M=M_\lambda$, and $L=mf^*A$.

\begin{proposition}\label{prop:picard-nine}
For every such smooth double cover, we have 
\[
h^1(M,\cO_M)=h^2(M,\cO_M)=0,
\quad \text{and} \quad \rho(M)=9.
\]
\end{proposition}
\begin{proof} 
By the Kawamata--Viehweg vanishing theorem and Serre duality, we have 
\begin{equation}\label{eq:negative-vanishing}
H^j(X,\cO_X(-kL))=0
\end{equation} 
for $0\leq j<4$ and $k>0$, because $L$ is nef and big (see, for example, \cite[Theorem 3.1.7]{fujino-foundations} or \cite[Corollary 1.7]{fujino-multiplication}). 
The higher cohomology of $\cO_X$ vanishes because $X$ is a complete
toric variety. Since $\pi_*\cO_M=\cO_X\oplus\cO_X(-L)$, we obtain $h^1(M, \mathcal O_M)=h^2(M, \mathcal O_M)=0$. The exponential sequence and Hodge decomposition
then show that $\rho(M)=h^{1,1}(M)$. 

The eigensheaf decomposition of the sheaf of differentials for a smooth double cover is
\[
\pi_*\Omega_M^1=\Omega_X^1\oplus
\bigl(\Omega_X^1(\log B_X)\otimes\cO_X(-L)\bigr), 
\]
where $B_X=B_T$ or $B_V$. 
We consider the short exact sequences:
\[
0\longrightarrow\Omega_X^1(-L)
\longrightarrow\Omega_X^1(\log B_X)(-L)
\longrightarrow\cO_{B_X}(-L)\longrightarrow0, 
\] 
and 
\[
0\longrightarrow\cO_X(-3L)\longrightarrow\cO_X(-L)
\longrightarrow\cO_{B_X}(-L)\longrightarrow0.
\]
By \eqref{eq:negative-vanishing}, we obtain 
$H^0(B_X, \cO_{B_X}(-L))=H^1(B_X, \cO_{B_X}(-L))=0$. Therefore, we have 
\[
H^1(X, \Omega_X^1(\log B_X)(-L))\simeq H^1(X, \Omega_X^1(-L)).
\]
We consider the twisted toric Euler sequence: 
\[
0\longrightarrow\Omega_X^1(-L)
\longrightarrow\bigoplus_{i=1}^{12}\cO_X(-D_i-L)
\longrightarrow\cO_X(-L)^{\oplus8}\longrightarrow0.
\]
Using \eqref{eq:negative-vanishing} once more, and then the sequences
\[
0\longrightarrow\cO_X(-D_i-L)\longrightarrow\cO_X(-L)
\longrightarrow\cO_{D_i}(-L)\longrightarrow0,
\]
we get
\begin{equation}\label{eq:picard-formula}
\rho(M)=\rho(X)+\sum_{i=1}^{12}
h^0(D_i,\cO_{D_i}(-L)).
\end{equation}
Since $L$ is nef, each restriction $L|_{D_i}$ is also nef. If $\cO_{D_i}(-L)$ has a nonzero section, then $L|_{D_i}$ is trivial.
Conversely, if $L|_{D_i}$ is trivial, it contributes one to the sum in \eqref{eq:picard-formula}.
Since $A$ is ample, $L|_{D_i}$ is trivial precisely when $f(D_i)$ is a
point. Lemma~\ref{lem:divisor-images} says that this happens only for
$i=5$. Therefore, $\rho(M)=8+1=9$.
\end{proof}

\subsection{The extra divisor class}
The branch divisor $B$ avoids $p=f(D_{5})$. Therefore, we have 
\[
\pi^*D_{5}=E_++E_-,\quad
E_+\cap E_-=\varnothing,\quad E_\pm\simeq D_{5}.
\]
Let $\iota$ be the deck involution and put $\Delta=E_+-E_-$.
Then $\iota^*\Delta=-\Delta$.

The wall $\langle u_{1},u_{2},u_{5}\rangle$ in the fan of either smooth endpoint corresponds to the relation 
$u_{2}+u_{3}+u_{4}=u_{1}+2u_{5}$. 
The corresponding invariant curve $C_0$ lies in $D_{5}$ and satisfies
$D_{5}\cdot C_0=-2$. For its copy $C_+\subset E_+$, this gives 
$\Delta\cdot C_+=-2$. 
In particular, $\Delta$ is numerically nonzero. For a finite surjective morphism, the pullback map on numerical divisor classes is injective, and its image is invariant under $\iota$. Hence, by Proposition~\ref{prop:picard-nine}, we obtain 
\[
N^1(M)_{\Q}=\pi^*N^1(X)_{\Q}\oplus\Q[\Delta].
\]
This is a splitting of vector spaces over $\Q$, not a claim about an integral direct-sum decomposition of Picard groups. Indeed, $2[E_+]=\pi^*[D_{5}]+[\Delta]$.

\section{A four-step elementary flop sequence}\label{sec:four-step}
The existence of some finite flop factorization follows from
\cite[Theorem~1]{Kawamata2008Flops}, since $M_T$ and $M_V$ are birational
smooth projective minimal models. We give a more precise factorization
by separating the two sheets of each designated toric flop.

\subsection{The two exceptional planes and their curve classes}
The exceptional locus of the flopping contraction $X_T\to P_{TS}$ 
is 
$V(\langle u_{3},u_{8}\rangle)\simeq\mathbb P^2$. 
Similarly, we see that 
the exceptional locus of the flopping contraction $X_V\to P_{SV}$ is
$V(\langle u_{6},u_{7}\rangle)\simeq\mathbb P^2$. Both planes map to $p\in W$.
These descriptions follow either from the stars of the two cones or from
the extremal relations $\beta_3,\gamma_3$. In each case, $D_{5}$ restricts
to $\cO_{\mathbb P^2}(1)$. Its intersection with a line is the coefficient
$1$ of $u_{5}$ in the primitive relation.
For $\lambda=T,V$, choose an ample $\Q$-Cartier divisor $H_\lambda$
whose intersection with a line in the corresponding exceptional plane is
one. Such a divisor is obtained by rescaling any ample Cartier divisor.

Since the branch divisor avoids $p$, the inverse image of either plane is the
disjoint union of two planes $\mathbb P^2_+,\mathbb P^2_-$. For their
line classes $\ell_+,\ell_-$ and $H=H_T$ or $H_V$, we have 
the following: 
\begin{equation}\label{eq:line-intersections}
\begin{array}{|c|rr|}
\hline
&\ell_+&\ell_-\\\hline
\pi^*H&1&1\\
\Delta&1&-1\\\hline
\end{array}
\end{equation}
Here, the two sheets are labelled consistently with $E_+,E_-$. The
intersection with $D_{5}$ takes place on the corresponding sheet.
The two curve classes are independent.

The base-changed contraction
\[
q_T:M_T\longrightarrow\widetilde P_{TS}:=P_{TS}\times_W\wtW
\]
contracts precisely these two planes. Its relative Mori cone is generated by
$\ell_+,\ell_-$. Therefore, we have 
$\rho(M_T/\widetilde P_{TS})=2$. 
The same statement holds for $M_V\to\widetilde P_{SV}$.
The targets are projective and normal. 
Note that the two original flops do not remain elementary under the base change.

\subsection{Separating the two rays}
Over $\widetilde P_{TS}$, the divisors
$\pi^*H_T+\Delta$ and $\pi^*H_T-\Delta$ 
are nef, with degrees $(2,0)$ and $(0,2)$ on the two rays. They are
relatively semiample by the relative basepoint-free theorem because 
$K_{M_T}$ is numerically trivial over $\widetilde P_{TS}$. 
A sufficiently divisible positive multiple of each divisor therefore defines 
a projective small contraction of exactly one of the planes. Each of these contractions admits its flop by 
\cite{BirkarEtAl2010Minimal}. 

Even after this operation, the two exceptional loci over
$\widetilde P_{TS}$ remain disjoint. The first has been flopped, and
the other remains a plane. The remaining plane gives the other
extremal ray on which the transform of an ample divisor from the
simultaneously flopped model is negative. The first flopped model is
$\Q$-factorial and terminal. We can apply similar arguments again. 
By flopping that ray, we obtain 
$M_S:=X_S\times_W\wtW$. 
By the same construction starting from $M_V$ and using $\gamma_3$, we can 
connect $M_V$ to $M_S$ by two elementary flops. By reversing this second pair, 
we finally obtain 
\begin{equation}\label{eq:four-step}
M_T\dashrightarrow M_1\dashrightarrow M_S
\dashrightarrow M_3\dashrightarrow M_V.
\end{equation}
We note that all contraction targets and all models in this sequence are projective.

\subsection{Singularities along the sequence}
The unique singular cone of $X_S$ lies over $p$, where the double cover is
unramified. Therefore, it produces two singularities of type
$\frac12(1,1,1,1)$ on $M_S$. There are no
other singularities on $M_S$.

Each separated flop changes only one of the two disjoint neighborhoods.
Hence, the numbers of singularities of this type in
\eqref{eq:four-step} are 
$0$, $1$, $2$, $1$, and $0$, respectively. 
The intermediate models are $\Q$-factorial and terminal. 
Their canonical
divisors are still pulled back from $K_{\wtW}$. Therefore, they remain nef.
This proves Theorem~\ref{thm:main} (ii) without
asserting that four steps are the shortest possible length.

\section{Local analytic Picard groups and toric modifications}\label{sec:local-picard}
To exclude arbitrary smooth projective factorizations, we need a local comparison of line bundles and models. We
state the comparison for a general full-dimensional strongly convex
rational cone $\sigma$.

\subsection{Line bundles near the closed fiber}
Let $g:Z_\Theta\to U_\sigma$ be a smooth projective toric resolution,
let $p$ be the closed orbit of $U_\sigma$, and let
$F=g^{-1}(p)_{\red}$.

\begin{lemma}\label{lem:local-picard}
Restriction induces a surjection
\[
\Pic(Z_\Theta)\longrightarrow
\varinjlim_{\Omega\ni p}\Pic(g^{-1}(\Omega)),
\]
where $\Omega$ runs through analytic neighborhoods of $p$ and the Picard
groups on the right are analytic. In particular, every line bundle on
the germ along $F$ is represented by the restriction of a torus-invariant
Cartier divisor on $Z_\Theta$.
\end{lemma}
\begin{proof}
We note that $U_\sigma$ has only rational singularities since it is toric. 
Therefore, if $\Omega$ is a sufficiently small Stein open neighborhood of $p$, 
we have 
\[
H^1(g^{-1}(\Omega),\cO)=H^2(g^{-1}(\Omega),\cO)=0.
\]
Hence, the exponential sequence identifies
\begin{equation}\label{eq:local-c1}
\Pic(g^{-1}(\Omega))\xrightarrow{\ c_1\ }
H^2(g^{-1}(\Omega),\Z)
\end{equation}
isomorphically. We take the direct limit and apply topological
proper base change. Then we have 
\begin{equation}\label{eq:proper-base-change}
\varinjlim_{\Omega\ni p}H^2(g^{-1}(\Omega),\Z)
\simeq (R^2g_*\Z)_p\simeq H^2(F,\Z).
\end{equation}

The variety $Z_\Theta$ is smooth and semiprojective (see, for example, 
\cite[Proposition 7.2.9]{CoxLittleSchenck2011Toric}). Its core, the union
of its complete invariant orbit closures, is precisely $F$. Indeed,
a complete invariant orbit closure must map to the unique closed orbit
of the affine base, and every orbit closure in $F$ is complete.
By \cite[Theorems~3.2 and~3.5]{HauselSturmfels2002Toric}, the inclusion
$F\hookrightarrow Z_\Theta$ induces an isomorphism on integral cohomology.
Moreover, $H^2(Z_\Theta,\Z)$ is generated by the first Chern classes of
invariant divisors by \cite[Proposition~2.11]{HauselSturmfels2002Toric}.

We represent the Chern class of a given element 
of \[\varinjlim_{\Omega\ni p}\Pic(g^{-1}(\Omega))\] in $H^2(F,\Z)$
by an integral combination of those invariant classes.
Equation~\eqref{eq:proper-base-change} implies that the two classes agree
on a sufficiently small neighborhood, and injectivity of
\eqref{eq:local-c1} then identifies the line bundles there. This proves
the asserted surjectivity.
\end{proof}

It is the comparison with the \emph{analytic germ}, not just the
algebraic Picard group of $Z_\Theta$, that is needed here. The next
lemma is deduced from this comparison. 

\subsection{Recovering a small modification by relative Proj}
The following lemmas will play an important role in the proof of 
Theorem~\ref{thm:main}. 

\begin{lemma}\label{lem:toric-modification}
Consider a commutative diagram
\[
\begin{tikzcd}
\mathcal X \arrow[rr, dashed] \arrow[dr] & & \mathcal X^{+} \arrow[dl] \\
& \mathcal Y &
\end{tikzcd}
\]
of normal algebraic varieties, where $\mathcal X\to \mathcal Y$ and $\mathcal X^{+}\to \mathcal Y$ are projective birational morphisms, and $\mathcal X\dashrightarrow \mathcal X^{+}$ is an isomorphism in codimension one. Suppose that the germ of the morphism $\mathcal X\to \mathcal Y$ over a point $y\in \mathcal Y$ is analytically isomorphic to the germ of a projective toric model $(Z_{\Theta}/U_{\sigma})_{p}$. Then the germ of $\mathcal X^{+}\to \mathcal Y$ over $y$ is analytically isomorphic to the germ of a projective toric model $(Z^{+}/U_{\sigma})_{p}$. The fan of $Z^{+}$ has the same rays as $\Theta$. Furthermore, if $\mathcal X^{+}$ is smooth, then the toric model $Z^{+}$ is also smooth.
\end{lemma}

\begin{proof}
Let $D^{+}$ be a Cartier divisor on $\mathcal X^{+}$ that is ample over $\mathcal Y$. Let $D$ be its strict transform on $\mathcal X$. 
Note that $\mathcal X$ is smooth over a small analytic open neighborhood of $y$. 
Hence, by Lemma~\ref{lem:local-picard}, the line bundle $\mathcal{O}_{\mathcal X}(D)$ is represented by an invariant Cartier divisor $D_{0}=\sum_{\rho}d_{\rho}D_{\rho}$ on $Z_{\Theta}$ over a sufficiently small analytic open neighborhood of $y$. 
We fix one such isomorphism and use its tensor powers for all degrees. 
Since $\mathcal X\dashrightarrow \mathcal X^{+}$ is an isomorphism in codimension one, the relative section algebra of $D^{+}$ over the germ coincides with the analytification of
\[
\mathcal{R}(D_{0})=\bigoplus_{r\ge 0}g_{*}\mathcal{O}_{Z_{\Theta}}(rD_{0}).
\]
It is well known that this algebra is finitely generated (see, for example, \cite[Theorem 4.7]{toric-mori-theory}).
Its relative $\operatorname{Proj}$ gives a normal projective toric model $Z^{+}$ over $U_{\sigma}$. This completes the proof of the first assertion.

Every torus invariant point of a projective toric birational model over $U_{\sigma}$ maps to $p$. 
Since every invariant prime divisor contains such a torus invariant point, any added or removed ray would be detected as a divisorial change in the germ. 
The smallness assumption on $\mathcal X \dashrightarrow \mathcal X^+$ rules this out. 
This means that the rays of $Z^+$ and $\Theta$ coincide. 
Finally, the smoothness of the germ tests the nonsingularity of every maximal cone. 
Thus, if $\mathcal X^+$ is smooth, $Z^+$ is also smooth.
\end{proof}

The same argument also applies to contractions.

\begin{lemma}\label{lem:toric-contraction}
Consider a commutative diagram
\[
\begin{tikzcd}
\mathcal X \arrow[rr, "\varphi"] \arrow[dr] & & \mathcal X' \arrow[dl] \\ 
& \mathcal Y &
\end{tikzcd}
\]
of normal algebraic varieties, where $\mathcal X \to \mathcal Y$ and $\mathcal X^{\prime} \to \mathcal Y$ are projective birational morphisms. Suppose that the germ of the morphism $\mathcal X \to \mathcal Y$ over a point $y \in \mathcal Y$ is analytically isomorphic to the germ of a projective toric model $(Z_{\Theta}/U_{\sigma})_{p}$. Then the germ of the morphism $\varphi \colon \mathcal X \to \mathcal X^{\prime}$ over $y$ is analytically isomorphic to the germ of a projective toric contraction over $U_{\sigma}$.
\end{lemma}

\begin{proof}
Let $H$ be a Cartier divisor on $\mathcal X^{\prime}$ that is ample over $\mathcal Y$. 
By Lemma~\ref{lem:local-picard}, its pullback on $\mathcal X$ is represented by an invariant Cartier divisor $E$ on $Z_{\Theta}$. 
The degree of $E$ is nonnegative on every complete invariant curve in the closed fiber. 
Every complete invariant curve maps to the unique torus invariant point $p$ of the affine base $U_\sigma$. Since $|\Theta|=\sigma$ is full-dimensional and convex, these curves generate $\NE(Z_\Theta/U_\sigma)$ by \cite[Theorem~6.3.20 (b) and Corollary~6.3.21]{CoxLittleSchenck2011Toric}. 
Thus, $E$ is relatively nef. 
In the toric setting, it is therefore relatively semiample by \cite[Theorem~6.3.12]{CoxLittleSchenck2011Toric}. 
The relative $\operatorname{Proj}$ of its section algebra analytically recovers the germ of $\mathcal X^{\prime} \to \mathcal Y$. 
This completes the proof.
\end{proof}

\section{Excluding smooth factorizations over the canonical model}\label{sec:obstruction}
We first prove the assertion of Theorem~\ref{thm:main} with the canonical model marked.

\subsection{All flops occur over the same base}
Suppose there were a sequence
\begin{equation}\label{eq:hypothetical-chain}
M_T=M_0\dashrightarrow M_1\dashrightarrow\cdots
\dashrightarrow M_s=M_V
\end{equation}
of elementary flops through smooth projective varieties. 
Note that each $K_{M_i}$ is
nef and big. Hence it is semiample by the basepoint-free theorem. 
Of course, the canonical ring is unchanged. Therefore, every $M_i$ has a canonical
morphism to $\wtW$, and its canonical divisor is the pullback of the
same ample $K_{\wtW}$. 
Let $q_i:M_i\to P_i$ be an elementary flopping contraction in the sequence.
It is not difficult to see that the entire sequence is
over $\wtW$.

\subsection{A split analytic neighborhood}
Let 
$\Omega$ be a small analytic open neighborhood of $p$. 
Then we have 
\[
\pi^{-1}(\Omega)=\Omega_+\sqcup\Omega_-,\quad
\Omega_\pm\simeq\Omega 
\]
since $p\notin B$. 
This is an analytic splitting. We write
$\widetilde p_\pm$ for the two points over $p$.
On either sheet, the endpoint germs are precisely
$(Z_T/U_\sigma)_p$ and $(Z_V/U_\sigma)_p$.

\subsection{Numerical localization and the elementary condition}\label{subsec:numerical-localization}
A global elementary contraction need not remain elementary on an
arbitrary analytic neighborhood. We now establish the necessary
numerical control in this example.

\begin{lemma}\label{lem:numerical-localization}
Suppose a smooth projective variety $M_i$ over $\wtW$ has, on one sheet,
a toric germ $(Z_\Theta/U_\sigma)_p$, and all invariant Cartier divisors
of this local model are restrictions of global Cartier divisors on
$M_i$. Then the restriction map
\[
N^1(M_i)_{\Q}\longrightarrow N^1(Z_\Theta/U_\sigma)_{\Q}
\]
is well defined and surjective. Its dual identifies the local curve
space injectively with a subspace of $N_1(M_i)_{\Q}$.
\end{lemma}
\begin{proof} 
We can restrict 
any global line bundle to a line bundle on the analytic germ. 
By Lemma~\ref{lem:local-picard}, it is represented by an algebraic invariant
divisor on $Z_\Theta$. Two such representatives agree numerically on
the torus invariant complete curves in $F=g^{-1}(p)_{\red}$. These curves generate the relative curve
space (see, for example, \cite[Corollary 6.3.21]{CoxLittleSchenck2011Toric}). Therefore, the resulting relative numerical class is well defined. 
We note that a numerically trivial global divisor restricts to the zero
relative class. This constructs the desired map. 
It is surjective by the extension hypothesis. Its dual is consequently
injective. 
\end{proof}

For the initial model $M_T$, the extension hypothesis is satisfied by
\begin{equation}\label{eq:global-invariant-divisors}
G_{T,j}=\pi_T^*D_j\quad(1\leq j\leq12).
\end{equation}
On one sheet, $G_{T,j}$ restricts, for $1\leq j\leq8$, to the invariant
divisor for $u_j$ on $Z_T$. We do not require the irreducibility of $G_{T,j}$.
Strict transforms of these divisors remain Cartier on any smooth
intermediate model.

We consider the first step of \eqref{eq:hypothetical-chain} that is
nontrivial on this marked germ. Before that step the germ is still
$Z_T$, and the strict transforms of \eqref{eq:global-invariant-divisors}
are still its invariant divisors locally. Lemmas~\ref{lem:toric-contraction}
and~\ref{lem:toric-modification} turn this step into a toric contraction
and its toric small modification over $U_\sigma$.

Let $F_i$ be the contracted face of $\NE(Z_T/U_\sigma)$.
Every curve in this face maps into the one-dimensional cone
$\NE(M_i/P_i)$. By Lemma~\ref{lem:numerical-localization}, the map on
curve spaces is injective. Hence we have $\dim F_i\leq1$. Since the step is
nontrivial, it is an elementary toric flop locally. 
However, \eqref{eq:relative-cones} and the criterion recalled in Section~\ref{completemori3.1} 
show that
each of the four possible local extremal flops 
produces a singular variety. 
This contradicts the smoothness of $M_{i+1}$.

Such a first nontrivial step must exist because the final germ $(Z_V/U_\sigma)_p$ is distinct from 
$(Z_T/U_\sigma)_p$. Therefore, we have proved the following.
\begin{proposition}\label{prop:marked-nonconnectivity}
For every branch satisfying the basic branch conditions, the specified
birational map $M_T\dashrightarrow M_V$ over $\wtW$ has no factorization
into elementary flops through smooth projective varieties.
\end{proposition}

\section{Removing the marking}\label{sec:automorphisms}
We now choose the branch divisor $B$ so that an alternative isomorphism of the canonical model cannot bypass Proposition~\ref{prop:marked-nonconnectivity}.

\begin{lemma}\label{lem:branch-stabilizer}
For all sufficiently large and divisible integers $m$, a nonempty
Zariski open subset of $|2mA|$ consists of divisors $B$ with
$\Aut(W,B)=\{1\}$.
\end{lemma}

\begin{proof}
The pullback map $\Pic(W)\to\Pic(X_T)$ is injective since 
$(f_T)_*\cO_{X_T}=\cO_W$. Since $\Pic(X_T)\simeq\Z^8$, the group $\Pic(W)$ is torsion-free. 
Therefore, an automorphism preserving $B\sim 2mA$ preserves the class of $A$. 
Hence, we have $\Aut (W, B)\subset \Aut (W, \mathcal O_W(A))$. 

We set $G:=\Aut(W,\mathcal O_W(A))$. 
It is well known that $G$ is a finite-type algebraic subgroup of a projective linear group. Let $d = \dim G$ and fix an integer $q>d$. 
If $m$ is sufficiently large and divisible, the sections of 
$2mA$ separate arbitrary $2q$ distinct points independently. 

For any nonidentity element $g\in G$, we can choose $q$ disjoint pairs of distinct points in the moved locus:
\[
x_1,gx_1,\ldots,x_q,gx_q.
\]
We set $V_m=H^0(W,2mA)$ and lift the projective action of $g$ to a linear automorphism of $V_m$. If a section $s \in V_m$ represents a fixed point in $\mathbb{P}(V_m)$, it must be an eigenvector of $g$. For any fixed eigenvalue, the eigenvector condition imposes one independent linear relation on the values of the section at each pair $\{x_i, gx_i\}$. Since the sections can separate these $2q$ points independently, each eigenspace has codimension at least $q$ in $V_m$. Because there are only finitely many eigenvalues, the fixed-point locus of $g$ in $\mathbb P(V_m)$ also has codimension at least $q$.

Consider the algebraic incidence variety
\[
I=\{(g,[s])\in(G\setminus\{1\})\times\mathbb P(V_m)\mid g^*[s]=[s]\}.
\]
By the above argument, the fibers of the projection $I \to G\setminus\{1\}$ have dimension at most $\dim\mathbb P(V_m)-q$. Thus, the total dimension of $I$ is bounded by
\[
\dim I \leq \dim G + \dim\mathbb P(V_m) - q = d + \dim\mathbb P(V_m) - q.
\]
Since we chose $q > d$, we strictly have $\dim I < \dim\mathbb P(V_m)$. Consequently, the closure of the projection of $I$ onto $\mathbb P(V_m)$ is a proper closed subset. 

For any section $[s]$ outside this closed subset (which forms a nonempty Zariski open set), its projective stabilizer in $G$ is trivial. Finally, on the normal complete variety $W$, two nonzero sections defining the same zero divisor differ by a nonzero constant multiplier, meaning they correspond to the same point in $\mathbb P(V_m)$. Therefore, the corresponding general zero divisor $B \in |2mA|$ has no nonidentity automorphism.
\end{proof}

This open condition can be imposed together with the basic branch conditions. We also keep $K_W+mA$ very ample, as chosen above.

\begin{lemma}\label{lem:cover-automorphisms}
For such a branch divisor $B$,
\[
\Aut(\wtW)=\{1,\iota\},
\]
where $\iota$ is the deck involution.
\end{lemma}
\begin{proof}
We first note that 
\[
\pi_*\mathcal O_{\wtW}(K_{\wtW})=\mathcal O_W(K_W)\oplus \mathcal O_W(K_W)\otimes \mathcal L 
\] 
by construction. 
Since $W$ is a complete toric variety, we have $H^0(W,K_W)=0$. Therefore, 
we obtain  
\[
\begin{aligned}
H^0(\wtW,K_{\wtW})
&=H^0(W,K_W+mA)\oplus H^0(W,K_W)\\
&=H^0(W,K_W+mA).
\end{aligned}
\]
Thus, the morphism defined by the canonical linear system on $\wtW$ factors as 
\[
\wtW\xrightarrow{\pi}W\hookrightarrow
\mathbb P (H^0(W,K_W+mA)^*).
\]
An automorphism of $\wtW$ must preserve this morphism and its image. 
Therefore, it descends to an automorphism of $W$ preserving $B$. 
By Lemma~\ref{lem:branch-stabilizer}, this descended automorphism is the identity. 
Because $\wtW \to W$ is a double cover, its only automorphisms over $W$ are the identity and the deck involution. This proves the claim.
\end{proof}

Let us prove Theorem~\ref{thm:main}. 

\begin{proof}[Proof of Theorem~\ref{thm:main}]
Assertion (i) follows from the construction and Proposition~\ref{prop:picard-nine}.
Assertion (ii) was proved in Section~\ref{sec:four-step}. Suppose that there is a sequence of smooth projective elementary flops from $M_T$ to a variety isomorphic to $M_V$.
This sequence induces an identification of the canonical rings.
By comparing this with the specified canonical morphism of $M_V$, we obtain an automorphism of $\wtW$.
By Lemma~\ref{lem:cover-automorphisms}, this automorphism is the identity or $\iota$.
If it is the identity, this contradicts Proposition~\ref{prop:marked-nonconnectivity}. On the other hand, the deck involution lifts to $M_V$ by the fiber-product construction.
If we compose the terminal isomorphism with this lift, the identification of the canonical model becomes the identity.
This gives the same contradiction.
This proves (iii).
\end{proof}

\begin{remark}\label{rem:scope-conclusion}
This construction provides families of examples, rather than a single branch polynomial with specific coefficients. 
If we fix a sufficiently large $m$, the conditions used above define a nonempty Zariski open subset. 
The coefficient $4$ in $K_W+4A$ just guarantees the ampleness. 
It does not mean that $m=4$ already satisfies all the conditions for the branch divisor and the automorphisms. 
Furthermore, we do not claim that the Picard numbers are minimal, nor that the four-step factorization is optimal.
\end{remark}

\appendix
\section{Exact decompositions for the relative Mori cones}\label{app:local-cones}
The fans of $Z_T$ and $Z_V$ are 
$\Sigma_T|_\sigma$ and $\Sigma_V|_\sigma$, respectively. 
This appendix makes the local extremal-ray argument independent of an unlisted cone computation. For each model, the table contains all nine distinct primitive integral relation vectors of the interior walls. The last column shows the number of walls that give that vector. In both cases, the sum is nineteen. The middle column shows the coordinates of the vector. We use the ordered basis displayed above the table. Every coordinate is nonnegative. Also, every basis vector appears. The basis vectors are linearly independent. This proves \eqref{eq:relative-cones}.

\subsection*{The model $Z_T$}
Ordered basis: $(\alpha_2,\alpha_4,\beta_2,\beta_3)$. 
\begin{center}
\small
\begin{tabular}{|l|c|c|}\hline
Relation & Coordinates & Walls\\\hline
$u_{2}+u_{3}+u_{4}=u_{1}+2u_{5}$ & $(0,1,0,0)$ & 3\\
$u_{4}+u_{5}+u_{7}=2u_{3}+u_{8}$ & $(0,0,0,1)$ & 3\\
$u_{3}+u_{6}=u_{4}+u_{5}$ & $(1,1,0,0)$ & 1\\
$u_{2}+u_{6}+u_{7}=3u_{5}$ & $(2,2,1,1)$ & 1\\
$u_{2}+u_{4}+u_{8}=2u_{5}+u_{6}$ & $(0,0,1,0)$ & 3\\
$u_{2}+u_{3}+u_{8}=3u_{5}$ & $(1,1,1,0)$ & 3\\
$u_{1}+u_{5}+u_{6}=u_{2}+2u_{4}$ & $(1,0,0,0)$ & 3\\
$u_{1}+u_{7}=2u_{3}$ & $(1,0,1,1)$ & 1\\
$u_{1}+u_{8}=u_{4}+u_{5}$ & $(1,0,1,0)$ & 1\\
\hline
\end{tabular}
\end{center}

\subsection*{The model $Z_V$}
Ordered basis: $(\alpha_2,\alpha_4,\gamma_2,\gamma_3)$. 
\begin{center}
\small
\begin{tabular}{|l|c|c|}\hline
Relation & Coordinates & Walls\\\hline
$u_{2}+u_{3}+u_{4}=u_{1}+2u_{5}$ & $(0,1,0,0)$ & 3\\
$u_{4}+u_{5}+u_{8}=2u_{6}+u_{7}$ & $(0,0,0,1)$ & 3\\
$u_{3}+u_{6}=u_{4}+u_{5}$ & $(1,1,0,0)$ & 1\\
$u_{2}+u_{4}+u_{7}=u_{3}+2u_{5}$ & $(0,1,1,0)$ & 1\\
$u_{2}+u_{6}+u_{7}=3u_{5}$ & $(1,2,1,0)$ & 4\\
$u_{2}+u_{4}+u_{8}=2u_{5}+u_{6}$ & $(1,2,1,1)$ & 1\\
$u_{2}+u_{3}+u_{8}=3u_{5}$ & $(2,3,1,1)$ & 1\\
$u_{1}+u_{5}+u_{6}=u_{2}+2u_{4}$ & $(1,0,0,0)$ & 3\\
$u_{1}+u_{7}=2u_{3}$ & $(0,0,1,0)$ & 2\\
\hline
\end{tabular}
\end{center}

The verifier \path{verification/relative.py} reconstructs the interior walls. It solves these coordinates over the rational numbers. It also checks the independence of the four generators. It reconstructs every flop discussed here. Finally, it checks the two new cone multiplicities.

\section{The common base and its Cartier data}\label{app:base}
For the coefficient calculations in this appendix, we use the auxiliary vector
\begin{equation}\label{eq:base-height}
a=(a_1,\ldots,a_{12})=(0,0,0,0,0,0,0,0,6,12,2,1).
\end{equation}
For rays of $\Delta_W$, $a_i$ is the coefficient of $D_i^W$ in $A$; we set $a_3=a_5=0$ for the two inserted rays.

The table lists every maximal cone $\Gamma_\gamma$ explicitly. 
Put $I_\gamma=\{i\mid u_i\in\Gamma_\gamma\}$. The vectors $a_\gamma$ and $q_\gamma$ are the integral Cartier data of $A$ and $-K_W$, respectively, in the standard sign convention. Thus, we have
\[
\pair{a_\gamma}{u_i}=-a_i,\quad\pair{q_\gamma}{u_i}=-1
\quad(i\in I_\gamma),
\]
with $a$ as in \eqref{eq:base-height}. This information gives the full fan and both Cartier divisors.
\begin{center}
\small
\begin{tabular}{|r|l|l|l|}\hline
$\gamma$ & $\Gamma_\gamma$ & $a_\gamma$ & $q_\gamma$\\\hline
1 & $\langle u_{1},u_{2},u_{4},u_{6},u_{7},u_{8}\rangle$ & $(0,0,0,0)$ & $(0,0,0,-1)$\\
2 & $\langle u_{1},u_{2},u_{7},u_{10}\rangle$ & $(0,24,-12,12)$ & $(0,4,-2,1)$\\
3 & $\langle u_{1},u_{2},u_{6},u_{9},u_{10}\rangle$ & $(18,-12,6,12)$ & $(3,-2,1,1)$\\
4 & $\langle u_{1},u_{4},u_{7},u_{9}\rangle$ & $(0,0,6,6)$ & $(0,0,1,0)$\\
5 & $\langle u_{1},u_{7},u_{9},u_{10}\rangle$ & $(0,6,6,12)$ & $(0,1,1,1)$\\
6 & $\langle u_{1},u_{4},u_{6},u_{9}\rangle$ & $(6,-6,6,6)$ & $(1,-1,1,0)$\\
7 & $\langle u_{2},u_{6},u_{8},u_{10}\rangle$ & $(0,-12,-12,12)$ & $(0,-2,-2,1)$\\
8 & $\langle u_{2},u_{7},u_{8},u_{10}\rangle$ & $(-36,24,-48,12)$ & $(-6,4,-8,1)$\\
9 & $\langle u_{4},u_{6},u_{8},u_{9}\rangle$ & $(0,-6,6,6)$ & $(0,-1,1,0)$\\
10 & $\langle u_{4},u_{7},u_{8},u_{12}\rangle$ & $(-1,-1,2,2)$ & $(3,3,-6,-7)$\\
11 & $\langle u_{4},u_{7},u_{9},u_{11}\rangle$ & $(-2,-2,6,6)$ & $(0,0,1,0)$\\
12 & $\langle u_{4},u_{7},u_{11},u_{12}\rangle$ & $(-2,-2,5,5)$ & $(0,0,3,2)$\\
13 & $\langle u_{4},u_{8},u_{9},u_{11}\rangle$ & $(-2,-4,6,6)$ & $(0,-1,1,0)$\\
14 & $\langle u_{4},u_{8},u_{11},u_{12}\rangle$ & $(-2,-3,5,5)$ & $(0,-3,3,2)$\\
15 & $\langle u_{6},u_{8},u_{9},u_{10}\rangle$ & $(0,-12,6,12)$ & $(0,-2,1,1)$\\
16 & $\langle u_{7},u_{8},u_{10},u_{11}\rangle$ & $(-10,-2,4,12)$ & $(-2,0,0,1)$\\
17 & $\langle u_{7},u_{8},u_{11},u_{12}\rangle$ & $(-3,-2,4,5)$ & $(-3,0,0,2)$\\
18 & $\langle u_{7},u_{9},u_{10},u_{11}\rangle$ & $(-8,-2,6,12)$ & $(-1,0,1,1)$\\
19 & $\langle u_{8},u_{9},u_{10},u_{11}\rangle$ & $(-8,-4,6,12)$ & $(-1,-1,1,1)$\\
\hline
\end{tabular}
\end{center}
For the adjoint divisor, the Cartier data is $4a_\gamma-q_\gamma$. The verifier checks all inequalities
\[
\pair{4a_\gamma-q_\gamma}{u_j}>1-4a_j\quad(j\notin I_\gamma).
\]
It also checks all vertex and subdivision conditions used in Proposition~\ref{prop:common-base}.

\section{The six-step weak Fano path}\label{app:weak-fano}

Let $\bar X_j=X_{\bar\Sigma_j}$ be the model after $j$ flops in
\eqref{eq:weak-fano-path}, with
$\bar\Sigma_0=\bar\Sigma_T$ and $\bar\Sigma_6=\bar\Sigma_V$.
The primitive ray generators are $u_1,\ldots,u_{11}$ throughout.
The initial fan $\bar\Sigma_0$ is obtained from $\Sigma_T$ by replacing
its four maximal cones containing $u_{12}$ with
$\langle u_4,u_7,u_8,u_{11}\rangle$ and keeping all other maximal cones.

The following table describes the six flops explicitly.
At step $j$, we remove the four cones in the second column from
$\bar\Sigma_{j-1}$ and insert the four cones in the third column.
All other maximal cones remain unchanged, and we include all faces
of the resulting maximal cones to obtain $\bar\Sigma_j$.
The step numbers agree with those in \eqref{eq:weak-fano-path}.

\begingroup
\small
\setlength{\tabcolsep}{5pt}
\renewcommand{\arraystretch}{1.4}
\begin{longtable}{|c|l|l|}
\hline
Flop & Removed maximal cones & Inserted maximal cones\\
\hline
\endhead
\hline
\endfoot
(1) &
$\begin{gathered}
\langle u_3,u_4,u_7,u_8\rangle,
\langle u_3,u_4,u_7,u_9\rangle,\\
\langle u_4,u_7,u_8,u_{11}\rangle,
\langle u_4,u_7,u_9,u_{11}\rangle
\end{gathered}$
&
$\begin{gathered}
\langle u_3,u_4,u_8,u_{11}\rangle,
\langle u_3,u_4,u_9,u_{11}\rangle,\\
\langle u_3,u_7,u_8,u_{11}\rangle,
\langle u_3,u_7,u_9,u_{11}\rangle
\end{gathered}$\\
\hline
(2) &
$\begin{gathered}
\langle u_3,u_4,u_5,u_8\rangle,
\langle u_3,u_4,u_8,u_{11}\rangle,\\
\langle u_3,u_5,u_7,u_8\rangle,
\langle u_3,u_7,u_8,u_{11}\rangle
\end{gathered}$
&
$\begin{gathered}
\langle u_3,u_4,u_5,u_{11}\rangle,
\langle u_3,u_5,u_7,u_{11}\rangle,\\
\langle u_4,u_5,u_8,u_{11}\rangle,
\langle u_5,u_7,u_8,u_{11}\rangle
\end{gathered}$\\
\hline
(3) &
$\begin{gathered}
\langle u_3,u_4,u_5,u_{11}\rangle,
\langle u_3,u_4,u_9,u_{11}\rangle,\\
\langle u_3,u_5,u_7,u_{11}\rangle,
\langle u_3,u_7,u_9,u_{11}\rangle
\end{gathered}$
&
$\begin{gathered}
\langle u_3,u_4,u_5,u_7\rangle,
\langle u_3,u_4,u_7,u_9\rangle,\\
\langle u_4,u_5,u_7,u_{11}\rangle,
\langle u_4,u_7,u_9,u_{11}\rangle
\end{gathered}$\\
\hline
(4) &
$\begin{gathered}
\langle u_4,u_5,u_6,u_8\rangle,
\langle u_4,u_5,u_8,u_{11}\rangle,\\
\langle u_4,u_6,u_8,u_9\rangle,
\langle u_4,u_8,u_9,u_{11}\rangle
\end{gathered}$
&
$\begin{gathered}
\langle u_4,u_5,u_6,u_{11}\rangle,
\langle u_4,u_6,u_9,u_{11}\rangle,\\
\langle u_5,u_6,u_8,u_{11}\rangle,
\langle u_6,u_8,u_9,u_{11}\rangle
\end{gathered}$\\
\hline
(5) &
$\begin{gathered}
\langle u_4,u_5,u_6,u_{11}\rangle,
\langle u_4,u_5,u_7,u_{11}\rangle,\\
\langle u_5,u_6,u_8,u_{11}\rangle,
\langle u_5,u_7,u_8,u_{11}\rangle
\end{gathered}$
&
$\begin{gathered}
\langle u_4,u_5,u_6,u_7\rangle,
\langle u_4,u_6,u_7,u_{11}\rangle,\\
\langle u_5,u_6,u_7,u_8\rangle,
\langle u_6,u_7,u_8,u_{11}\rangle
\end{gathered}$\\[5pt]
\hline
(6) &
$\begin{gathered}
\langle u_4,u_6,u_7,u_{11}\rangle,
\langle u_4,u_6,u_9,u_{11}\rangle,\\
\langle u_6,u_7,u_8,u_{11}\rangle,
\langle u_6,u_8,u_9,u_{11}\rangle
\end{gathered}$
&
$\begin{gathered}
\langle u_4,u_6,u_7,u_8\rangle,
\langle u_4,u_6,u_8,u_9\rangle,\\
\langle u_4,u_7,u_8,u_{11}\rangle,
\langle u_4,u_8,u_9,u_{11}\rangle
\end{gathered}$\\
\end{longtable}
\endgroup

All seven fans have thirty-two nonsingular maximal cones.
The final fan is exactly $\bar\Sigma_V$.
The completeness and projectivity of the fans, the extremality of the
six relations in \eqref{eq:weak-fano-path}, and the nonnegative degrees
of all primitive relations at every stage are verified by exact
computer calculations; see Appendix~\ref{app:verification}.
Thus the table gives the six smooth projective elementary $K$-flops
in Proposition~\ref{prop:weak-fano-path}.

The primitive collections can also be computed successively using
\cite[Theorem~2.4 and Remark~2.5]{SatoTsuzuki2023AntiFlips}.
The resulting counts are the following: 
\begingroup
\setlength{\extrarowheight}{3pt}
\[
\begin{array}{|c|rrrrrrr|}
\hline
j&0&1&2&3&4&5&6\\\hline
\#\operatorname{PC}(\bar\Sigma_j)&26&25&22&23&22&25&26\\
\#\{\text{extremal primitive relations}\}&8&9&10&9&9&9&8\\
\hline
\end{array}
\]
\endgroup
In this table, $\operatorname{PC}(\bar\Sigma_j)$ stands for 
the set of primitive collections of $\bar\Sigma_j$. 
In each model the contractible primitive relations coincide with the
extremal ones. The final primitive collections agree with those of
$\bar\Sigma_V$ in Appendix~\ref{app:primitive-relations}.
The \texttt{primitive} section of \path{verification/data.json} contains the
full intermediate lists and exact extremality certificates. 
The script \path{verification/primitive.py} checks 
the changes in the primitive collections predicted by 
\cite[Theorem~2.4 and Remark~2.5]{SatoTsuzuki2023AntiFlips} 
against an independent enumeration from the maximal cones. 
The fan data and supporting heights for all models and contraction
targets are retained in the \texttt{weak\_fano} section of the same data file.

\section{Primitive relations and contractible classes}\label{app:primitive-relations}

We list all primitive relations of
$\Sigma_T,\Sigma_V,\bar\Sigma_T,\bar\Sigma_V$.
Exact computer calculations show that, for each of these four fans,
a primitive relation is contractible if and only if it is extremal.
We therefore divide each list into extremal primitive relations 
and the remaining primitive 
relations, which are neither contractible nor extremal.
In every relation, the generators on the left form the primitive collection. 
Contractibility is checked using
\cite[Proposition~3.4 and Theorem~2.2]{Casagrande2003Contractible},
and extremality is verified by the exact computations described in
Appendix~\ref{app:verification}.

\subsection{The fan $\Sigma_T$}\label{app:primitive-T}
There are 34 primitive relations.

\noindent\emph{Extremal primitive relations} (9):
\begingroup\small
\setlength{\jot}{2pt}
\[
\begin{gathered}
u_3+u_{12}=2u_4+2u_7+u_8,\quad 
u_9+u_{12}=u_4+3u_{11},\quad 
u_1+u_5+u_6=u_2+2u_4,\\ 
u_1+u_6+u_{10}=u_2+2u_9,\quad 
u_2+u_3+u_4=u_1+2u_5,\quad 
u_2+u_4+u_8=2u_5+u_6,\\ 
u_2+u_4+u_9=u_1+u_6,\quad 
u_4+u_5+u_7=2u_3+u_8,\quad 
u_4+u_7+u_8+u_{11}=u_{12}.
\end{gathered}
\]
\endgroup

\noindent\emph{Remaining primitive relations} (25):
\begingroup\small
\setlength{\jot}{2pt}
\[
\begin{gathered}
u_1+u_7=2u_3,\quad
u_1+u_8=u_4+u_5,\quad
u_1+u_{11}=u_3+u_4,\\
u_1+u_{12}=u_3+2u_4+u_7+u_8,\quad
u_2+u_{11}=2u_5,\quad
u_2+u_{12}=2u_3+u_5+2u_8,\\
u_3+u_6=u_4+u_5,\quad
u_3+u_{11}=u_4+u_7,\quad
u_4+u_{10}=u_9,\quad
u_5+u_9=u_4,\\
u_5+u_{10}=0,\quad
u_5+u_{11}=u_3+u_8,\quad
u_5+u_{12}=u_3+u_4+u_7+2u_8,\\
u_6+u_7=u_3+u_8,\quad
u_6+u_{11}=u_4+u_8,\quad
u_6+u_{12}=2u_4+u_7+2u_8,\\
u_{10}+u_{12}=3u_{11},\quad
u_2+u_3+u_8=3u_5,\quad
u_2+u_3+u_9=u_1+u_5,\\
u_2+u_4+u_7=u_3+2u_5,\quad
u_2+u_7+u_9=u_3+u_5,\quad
u_2+u_8+u_9=u_5+u_6,\\
u_3+u_8+u_9=u_4+u_{11},\quad
u_3+u_8+u_{10}=u_{11},\quad
u_7+u_8+u_9=2u_{11}.
\end{gathered}
\]
\endgroup

\subsection{The fan $\Sigma_V$}\label{app:primitive-V}
There are 34 primitive relations.

\noindent\emph{Extremal primitive relations} (9):
\begingroup\small
\setlength{\jot}{2pt}
\[
\begin{gathered}
u_1+u_7=2u_3,\quad
u_6+u_{12}=2u_4+u_7+2u_8,\quad
u_9+u_{12}=u_4+3u_{11},\\
u_1+u_5+u_6=u_2+2u_4,\quad
u_1+u_6+u_{10}=u_2+2u_9,\quad
u_2+u_3+u_4=u_1+2u_5,\\
u_2+u_4+u_9=u_1+u_6,\quad
u_4+u_5+u_8=2u_6+u_7,\quad
u_4+u_7+u_8+u_{11}=u_{12}.
\end{gathered}
\]
\endgroup

\noindent\emph{Remaining primitive relations} (25):
\begingroup\small
\setlength{\jot}{2pt}
\[
\begin{gathered}
u_1+u_8=u_4+u_5,\quad
u_1+u_{11}=u_3+u_4,\quad
u_1+u_{12}=2u_4+u_6+2u_7,\\
u_2+u_{11}=2u_5,\quad
u_2+u_{12}=u_5+2u_6+2u_7,\quad
u_3+u_6=u_4+u_5,\\
u_3+u_8=u_6+u_7,\quad
u_3+u_{11}=u_4+u_7,\quad
u_3+u_{12}=2u_4+2u_7+u_8,\\
u_4+u_{10}=u_9,\quad
u_5+u_9=u_4,\quad
u_5+u_{10}=0,\quad
u_5+u_{11}=u_6+u_7,\\
u_5+u_{12}=u_4+u_6+2u_7+u_8,\quad
u_6+u_{11}=u_4+u_8,\quad
u_{10}+u_{12}=3u_{11},\\
u_2+u_3+u_9=u_1+u_5,\quad
u_2+u_4+u_7=u_3+2u_5,\quad
u_2+u_4+u_8=2u_5+u_6,\\
u_2+u_6+u_7=3u_5,\quad
u_2+u_7+u_9=u_3+u_5,\quad
u_2+u_8+u_9=u_5+u_6,\\
u_6+u_7+u_9=u_4+u_{11},\quad
u_6+u_7+u_{10}=u_{11},\quad
u_7+u_8+u_9=2u_{11}.
\end{gathered}
\]
\endgroup

\subsection{The fan $\bar\Sigma_T$}\label{app:primitive-barT}
There are 26 primitive relations.

\noindent\emph{Extremal primitive relations} (8):
\begingroup\small
\setlength{\jot}{2pt}
\[
\begin{gathered}
u_3+u_{11}=u_4+u_7,\quad
u_1+u_5+u_6=u_2+2u_4,\quad
u_1+u_6+u_{10}=u_2+2u_9,\\
u_2+u_3+u_4=u_1+2u_5,\quad
u_2+u_4+u_8=2u_5+u_6,\quad
u_2+u_4+u_9=u_1+u_6,\\
u_4+u_5+u_7=2u_3+u_8,\quad
u_7+u_8+u_9=2u_{11}.
\end{gathered}
\]
\endgroup

\noindent\emph{Remaining primitive relations} (18):
\begingroup\small
\setlength{\jot}{2pt}
\[
\begin{gathered}
u_1+u_7=2u_3,\quad
u_1+u_8=u_4+u_5,\quad
u_1+u_{11}=u_3+u_4,\quad
u_2+u_{11}=2u_5,\\
u_3+u_6=u_4+u_5,\quad
u_4+u_{10}=u_9,\quad
u_5+u_9=u_4,\quad
u_5+u_{10}=0,\\
u_5+u_{11}=u_3+u_8,\quad
u_6+u_7=u_3+u_8,\quad
u_6+u_{11}=u_4+u_8,\quad
u_2+u_3+u_8=3u_5,\\
u_2+u_3+u_9=u_1+u_5,\quad
u_2+u_4+u_7=u_3+2u_5,\quad
u_2+u_7+u_9=u_3+u_5,\\
u_2+u_8+u_9=u_5+u_6,\quad
u_3+u_8+u_9=u_4+u_{11},\quad
u_3+u_8+u_{10}=u_{11}.
\end{gathered}
\]
\endgroup

\subsection{The fan $\bar\Sigma_V$}\label{app:primitive-barV}
There are 26 primitive relations.

\noindent\emph{Extremal primitive relations} (8):
\begingroup\small
\setlength{\jot}{2pt}
\[
\begin{gathered}
u_1+u_7=2u_3,\quad
u_6+u_{11}=u_4+u_8,\quad
u_1+u_5+u_6=u_2+2u_4,\\
u_1+u_6+u_{10}=u_2+2u_9,\quad
u_2+u_3+u_4=u_1+2u_5,\quad
u_2+u_4+u_9=u_1+u_6,\\
u_4+u_5+u_8=2u_6+u_7,\quad
u_7+u_8+u_9=2u_{11}.
\end{gathered}
\]
\endgroup

\noindent\emph{Remaining primitive relations} (18):
\begingroup\small
\setlength{\jot}{2pt}
\[
\begin{gathered}
u_1+u_8=u_4+u_5,\quad
u_1+u_{11}=u_3+u_4,\quad
u_2+u_{11}=2u_5,\\
u_3+u_6=u_4+u_5,\quad
u_3+u_8=u_6+u_7,\quad
u_3+u_{11}=u_4+u_7,\\
u_4+u_{10}=u_9,\quad
u_5+u_9=u_4,\quad
u_5+u_{10}=0,\\
u_5+u_{11}=u_6+u_7,\quad
u_2+u_3+u_9=u_1+u_5,\quad
u_2+u_4+u_7=u_3+2u_5,\\
u_2+u_4+u_8=2u_5+u_6,\quad
u_2+u_6+u_7=3u_5,\quad
u_2+u_7+u_9=u_3+u_5,\\
u_2+u_8+u_9=u_5+u_6,\quad
u_6+u_7+u_9=u_4+u_{11},\quad
u_6+u_7+u_{10}=u_{11}.
\end{gathered}
\]
\endgroup

\section{Reproducibility and the scope of computation}\label{app:verification}
The verification code and exact computational data are provided in the
versioned supplement \cite{FujinoSato2026Verification}.
The ray identifiers in the code, data, and output agree with the
subscripts of $u_1,\ldots,u_{12}$ in this paper.
After extracting the archive, run the following command from its root:
\begin{center}
\verb|python3 verification/run_all.py|
\end{center}
The programs require Python~3.10 or later and only its standard library.
All mathematical checks use integers and exact rational arithmetic;
no floating-point tolerance is used.
The file \path{verification/data.json} contains the fixed fan data,
supporting heights, and exact certificates.
The checks cover the following assertions.
\begin{enumerate}
\item The fans of $X_T,X_S,X_V$, the two contraction targets,
and the common base $W$ are reconstructed from their generators.
The programs verify the defining supporting inequalities, polytope
boundedness and vertex lists, completeness, projectivity, cone
multiplicities, and the quotient type of the singular chart of $X_S$.
They also verify the common subdivisions through $W_1,W_2$,
their Picard numbers, and the Cartier, ampleness, and crepancy data
used in Proposition~\ref{prop:common-base}.
\item All walls and primitive collections are enumerated.
Nonnegative rational decompositions and exposing nef divisors certify
the complete Mori cones, rather than only a list of candidate rays.
The primitive relations in Appendix~\ref{app:primitive-relations}
are checked for contractibility using both criteria in
\cite[Proposition~3.4 and Theorem~2.2]{Casagrande2003Contractible}.
\item The relative Mori cones in Appendix~\ref{app:local-cones}
are reconstructed from their interior walls.
The programs check the local and global flopping rays, the
multiplicities of the new cones, the divisorial-ray counts,
and the obstruction to a smooth elementary toric flop from either
endpoint. They also check that the nonsingular charts of $X_S$
are covered by charts of the two smooth endpoints.
\item The finite toric inputs for the double-cover construction are
checked: the dimensions of the images of invariant divisors,
the relevant intersections with $D_5$, and the line classes of the
two exceptional planes used in \eqref{eq:line-intersections}.
\item The three-step six-ray path and its lattice identification with
\cite[Example~25]{Fujino2004Fequivalence} are checked.
For the six-step weak Fano path in Appendix~\ref{app:weak-fano},
the programs verify every model and contraction target, the
oriented extremal relations, and the changes in primitive collections
predicted by \cite[Theorem~2.4 and Remark~2.5]{SatoTsuzuki2023AntiFlips}.
\end{enumerate}

An optional transcription check compares the explicit numerical
displays in the manuscript with the same data:
\begin{center}
\verb|python3 verification/run_all.py --manuscript /path/to/main.tex|
\end{center}
The bibliography file should be in the same directory as the manuscript.
The script \path{verification/check_manuscript.py} compares the
ray coordinates, cone lists, primitive relations, relative wall tables,
Cartier data, and the two smooth paths. It also checks the cross-references
and citation keys. The \path{README.md} describes the commands, file
structure, and expected output.

The computer checks only the finite fan, lattice, and intersection inputs.
The analytic Picard lemma, the compatibility of relative section algebras,
the branch-divisor genericity argument, and the exclusion of non-toric
smooth factorizations are established by the geometric proofs in
Sections~\ref{sec:double-cover}--\ref{sec:automorphisms}.

\bibliographystyle{amsalpha}
\bibliography{ref}

@book{CoxLittleSchenck2011Toric,
  author    = {Cox, David A. and Little, John B. and Schenck, Henry K.},
  title     = {Toric Varieties},
  series    = {Graduate Studies in Mathematics},
  volume    = {124},
  publisher = {American Mathematical Society},
  address   = {Providence, RI},
  year      = {2011},
  MRNUMBER = {2810322},
  doi       = {10.1090/gsm/124}
}

@article{CoxVonRenesse2009Primitive,
  author  = {Cox, David A. and von Renesse, Christine},
  title   = {Primitive collections and toric varieties},
  journal = {Tohoku Math. J. (2)},
  volume  = {61},
  number  = {3},
  pages   = {309--332},
  year    = {2009},
  doi     = {10.2748/tmj/1255700197},
  eprint  = {0808.1836},
  MRNUMBER = {2568257},
  archivePrefix = {arXiv},
  primaryClass  = {math.AG}
}

@article{HauselSturmfels2002Toric,
  author  = {Hausel, Tam{\'a}s and Sturmfels, Bernd},
  title   = {Toric hyper{K}{\"a}hler varieties},
  journal = {Doc. Math.},
  volume  = {7},
  pages   = {495--534},
  year    = {2002},
  doi     = {10.4171/DM/130},
  MRNUMBER = {2015052},
  eprint  = {math/0203096},
  archivePrefix = {arXiv},
  primaryClass  = {math.AG}
}

@article {kawamata-d-equivalence,
    AUTHOR = {Kawamata, Yujiro},
     TITLE = {{$D$}-equivalence and {$K$}-equivalence},
   JOURNAL = {J. Differential Geom.},
  FJOURNAL = {Journal of Differential Geometry},
    VOLUME = {61},
      YEAR = {2002},
    NUMBER = {1},
     PAGES = {147--171},
      ISSN = {0022-040X,1945-743X},
   MRCLASS = {14F05 (14E05 14E30 18E30)},
  MRNUMBER = {1949787},
MRREVIEWER = {Ana\ Jerem\'{\i}as L\'{o}pez},
       URL = {http://projecteuclid.org.kyoto-u.idm.oclc.org/euclid.jdg/1090351323},
}

@article {orlov,
    AUTHOR = {Orlov, D. O.},
     TITLE = {Derived categories of coherent sheaves and equivalences
              between them},
   JOURNAL = {Uspekhi Mat. Nauk},
  FJOURNAL = {Uspekhi Matematicheskikh Nauk},
    VOLUME = {58},
      YEAR = {2003},
    NUMBER = {3(351)},
     PAGES = {89--172},
      ISSN = {0042-1316,2305-2872},
   MRCLASS = {14F05 (14A22 18E30)},
  MRNUMBER = {1998775},
MRREVIEWER = {Bal\'{a}zs\ Szendr\H{o}i},
       DOI = {10.1070/RM2003v058n03ABEH000629},
       URL = {https://doi-org.kyoto-u.idm.oclc.org/10.1070/RM2003v058n03ABEH000629},
}

@article{Kawamata2008Flops,
  author  = {Kawamata, Yujiro},
  title   = {Flops connect minimal models},
  journal = {Publ. Res. Inst. Math. Sci.},
  volume  = {44},
  number  = {2},
  pages   = {419--423},
  year    = {2008},
  MRNUMBER = {2426353},
  doi     = {10.2977/PRIMS/1210167332},
  eprint  = {0704.1013},
  archivePrefix = {arXiv},
  primaryClass  = {math.AG}
}

@article{BirkarEtAl2010Minimal,
  author  = {Birkar, Caucher and Cascini, Paolo and Hacon, Christopher D.
             and McKernan, James},
  title   = {Existence of minimal models for varieties of log general type},
  journal = {J. Amer. Math. Soc.},
  volume  = {23},
  number  = {2},
  pages   = {405--468},
  year    = {2010},
  MRNUMBER = {2601039},
  doi     = {10.1090/S0894-0347-09-00649-3},
  eprint  = {math/0610203},
  archivePrefix = {arXiv},
  primaryClass  = {math.AG}
}

@article {toric-mori-theory,
    AUTHOR = {Fujino, Osamu and Sato, Hiroshi},
     TITLE = {Introduction to the toric {M}ori theory},
   JOURNAL = {Michigan Math. J.},
  FJOURNAL = {Michigan Mathematical Journal},
    VOLUME = {52},
      YEAR = {2004},
    NUMBER = {3},
     PAGES = {649--665},
      ISSN = {0026-2285,1945-2365},
   MRCLASS = {14E30 (14M25)},
  MRNUMBER = {2097403},
MRREVIEWER = {James\ McKernan},
       DOI = {10.1307/mmj/1100623418},
       URL = {https://doi-org.kyoto-u.idm.oclc.org/10.1307/mmj/1100623418},
}

@article {fujino-multiplication,
    AUTHOR = {Fujino, Osamu},
     TITLE = {Multiplication maps and vanishing theorems for toric
              varieties},
   JOURNAL = {Math. Z.},
  FJOURNAL = {Mathematische Zeitschrift},
    VOLUME = {257},
      YEAR = {2007},
    NUMBER = {3},
     PAGES = {631--641},
      ISSN = {0025-5874,1432-1823},
   MRCLASS = {14M25 (14F10 14F17)},
  MRNUMBER = {2328817},
MRREVIEWER = {Jaros\l aw\ A.\ Wi\'{s}niewski},
       DOI = {10.1007/s00209-007-0140-5},
       URL = {https://doi-org.kyoto-u.idm.oclc.org/10.1007/s00209-007-0140-5},
}

@book {fujino-foundations,
    AUTHOR = {Fujino, Osamu},
     TITLE = {Foundations of the minimal model program},
    SERIES = {MSJ Memoirs},
    VOLUME = {35},
 PUBLISHER = {Mathematical Society of Japan, Tokyo},
      YEAR = {2017},
     PAGES = {xv+289},
      ISBN = {978-4-86497-045-7},
   MRCLASS = {14E30 (14F17)},
  MRNUMBER = {3643725},
MRREVIEWER = {Enrica\ Floris},
}

@article {fujino-fundamental,
    AUTHOR = {Fujino, Osamu},
     TITLE = {Fundamental theorems for the log minimal model program},
   JOURNAL = {Publ. Res. Inst. Math. Sci.},
  FJOURNAL = {Publications of the Research Institute for Mathematical
              Sciences},
    VOLUME = {47},
      YEAR = {2011},
    NUMBER = {3},
     PAGES = {727--789},
      ISSN = {0034-5318,1663-4926},
   MRCLASS = {14E30 (14C30 14F17)},
  MRNUMBER = {2832805},
MRREVIEWER = {Vladimir\ Lazi\'{c}},
       DOI = {10.2977/PRIMS/50},
       URL = {https://doi-org.kyoto-u.idm.oclc.org/10.2977/PRIMS/50},
}

@misc{Fujino2004Fequivalence,
  author       = {Fujino, Osamu},
  title        = {{$F$}-equivalence},
  howpublished = {Handwritten notes for the Kusatsu Algebraic Geometry Seminar},
  year         = {2004},
  note         = {June 29, 2004; \url{https://www.math.kyoto-u.ac.jp/~fujino/Ku-HP.pdf}},
  url          = {https://www.math.kyoto-u.ac.jp/~fujino/Ku-HP.pdf}
}

@article{Casagrande2003Contractible,
  author  = {Casagrande, Cinzia},
  title   = {Contractible classes in toric varieties},
  journal = {Math. Z.},
  volume  = {243},
  number  = {1},
  pages   = {99--126},
  year    = {2003},
  MRNUMBER = {1953051}, 
  doi     = {10.1007/s00209-002-0453-3},
  eprint  = {math/0111332},
  archivePrefix = {arXiv},
  primaryClass  = {math.AG}
}

@article{SatoTsuzuki2023AntiFlips,
  author  = {Sato, Hiroshi and Tsuzuki, Shigehito},
  title   = {Anti-flips of the blow-ups of the projective spaces at torus invariant points},
  journal = {Bull. Korean Math. Soc.},
  volume  = {60},
  number  = {6},
  pages   = {1705--1714},
  year    = {2023},
  doi     = {10.4134/BKMS.b220864},
  eprint  = {2212.08803},
  MRNUMBER = {4725865},
  archivePrefix = {arXiv},
  primaryClass = {math.AG}
}

@article {flops,
    AUTHOR = {Koll\'{a}r, J\'{a}nos},
     TITLE = {Flops},
   JOURNAL = {Nagoya Math. J.},
  FJOURNAL = {Nagoya Mathematical Journal},
    VOLUME = {113},
      YEAR = {1989},
     PAGES = {15--36},
      ISSN = {0027-7630,2152-6842},
   MRCLASS = {14E05 (14J30 32J25)},
  MRNUMBER = {986434},
MRREVIEWER = {Peter\ Nielsen},
       DOI = {10.1017/S0027763000001240},
       URL = {https://doi-org.kyoto-u.idm.oclc.org/10.1017/S0027763000001240},
}

@misc{FujinoSato2026Verification,
  author = {Fujino, Osamu and Sato, Hiroshi},
  title = {Verification code and data for {``Smooth minimal fourfolds not connected by smooth elementary flops''}},
  year = {2026},
  howpublished = {Zenodo, \url{https://doi.org/10.5281/zenodo.24340424}},
  note = {Version 1.0.0. Computational supplement},
}
\end{document}